\documentclass[journal]{IEEEtran}
\usepackage{cite}

\ifCLASSINFOpdf
  \usepackage[pdftex]{graphicx}
\else
\fi
\usepackage{amsmath}
\usepackage{amssymb}
\usepackage{mathtools}
\usepackage{setspace}
\usepackage{bm} 
\usepackage{units} 
\usepackage{soul}
\usepackage{balance}

\usepackage{amsthm} 

\newtheorem{proposition}{Proposition}

\theoremstyle{plain}
\newtheorem{remark}{Remark}

\usepackage{algorithm}
\usepackage[noend]{algpseudocode}
\usepackage[inline]{enumitem}
\usepackage{booktabs}
\usepackage{setspace}
\usepackage{supertabular}

\usepackage{pgfplots}
\pgfplotsset{compat=1.9}
\usepackage{tikzscale}

\usetikzlibrary{calc}

\usepackage{array}
\usepackage{url}
\usepackage[hidelinks]{hyperref}

\usepackage{xcolor}

\newcommand{\myDots}{\ifmmode\mathinner{\ldotp\kern-0.1em\ldotp\kern-0.1em\ldotp}\else.\kern-0.8em.\kern-0.8em.\fi}

\DeclareMathOperator*{\argmin}{\arg\min}

\newcommand*\circled[1]{\tikz[baseline=(char.base)]{\node[shape=circle,draw,inner sep=1pt] (char) {\footnotesize #1};}}

\newcommand{\eqlabel}[1]{\addtocounter{equation}{-1}\refstepcounter{equation}\label{#1}}

\newif\iflongversion

\longversiontrue

\begin{document}
%
\title{Strengthening the Group: Aggregated Frequency Reserve Bidding with ADMM}
%
%
%


\author{Felix~Rey,~\IEEEmembership{Member,~IEEE,}
        Xiaojing~Zhang,~\IEEEmembership{Member,~IEEE,}
        Sandro~Merkli,
        Valentina~Agliati,
        Maryam~Kamgarpour,~\IEEEmembership{Member,~IEEE,}
        and~John~Lygeros,~\IEEEmembership{Fellow,~IEEE}
\thanks{F.~Rey {\tt\small rey*}, S.~Merkli {\tt\small smerkli*}, M.~Kamgarpour {\tt\small mkamgar*} and J.~Lygeros {\tt\small lygeros*} are with the Automatic Control Laboratory at ETH Zurich, Switzerland, {\tt\small*@control.ee.ethz.ch}. Our research is supported by the European Research Council under ERC Starting Grant CONENE and by ABB Corporate Research under grant 2017-1224/01.”
X.~Zhang {\tt\small xiaojing.zhang@berkeley.edu} is with the MPC Lab at UC Berkeley, USA. V.~Agliati graduated from Politecnico di Milano in 2017 and contributed to this work during an exchange at ETH Zurich.\looseness=-1}}%

%
%

\markboth{Appears in IEEE Transactions on Smart Grid, DOI: 10.1109/TSG.2018.2841508, Print ISSN: 1949-3053,
Online ISSN: 1949-3061, \copyright 2018 IEEE}%
{Shell \MakeLowercase{\textit{et al.}}: Bare Demo of IEEEtran.cls for IEEE Journals}
%



\maketitle

\begin{abstract}
In a power grid, the electricity supply and demand must be balanced at all times to maintain the system's frequency.
In practice, the grid operator achieves this balance by procuring frequency reserves in an ahead-of-time market setting. During runtime, these reserves are then dispatched  whenever there is an imbalance in the grid. 
Recently, there has been an increasing interest in engaging electricity consumers, such  as plug-in electric vehicles or buildings, to offer such frequency reserves by exploiting their flexibility in power consumption. 
In this work, we focus on an aggregation of buildings that places a joint bid on a reserve market. 
The resulting shared decision is modeled as a large-scale optimization problem. 
Our main contribution is to show that the aggregation can  make its decision in a computationally efficient and conceptually meaningful way, by using the alternating direction method of multipliers (ADMM). 
The proposed approach exhibits several attractive features that include $(i)$~the computational burden is 
 distributed between the buildings; $(ii)$~the setup naturally provides privacy and flexibility; $(iii)$~the iterative algorithm can be stopped at any time, providing a feasible (though suboptimal) solution; and $(iv)$~the algorithm provides the foundation for a reward distribution scheme that strengthens the group.\looseness=-1
\end{abstract}


%
\IEEEpeerreviewmaketitle

\section{Introduction}
%
%
%
%
\IEEEPARstart{C}{ountries} around the world are integrating renewable energy sources into their power systems. 
The electricity production of weather-dependent renewable energy sources, such as wind turbines and photovoltaics, is volatile. 
This increases the uncertainty in the grid operation and challenges the balance of electricity supply and demand. 
 To maintain the active power balance and regulating the grid frequency, grid operators procure ancillary services, such as reserved generation capacities~\cite{rebours2007survey}. 
 Traditionally, these reserves are provided by conventional generators. 
However, the rapid integration of renewable energy sources calls for an increased balancing effort~\cite{makarov2009operational}. 
Recently, it has been shown that if properly aggregated and controlled, flexible loads can be engaged to provide frequency regulation at low cost and small environmental impact~\cite{callaway2011achieving}. 
Consumers with sizeable thermal storage, such as buildings, are particularly suited to contribute to the grid balance as their electricity consumption can be shifted in time with negligible impact on occupant \mbox{comfort~\cite{oldewurtel2013framework, Piette2012IntelligentBuilding}}. 

In many regions, consumers provide frequency reserves by offering them to the grid operator in a market setting. In Switzerland, the offered reserve size needs to be constant over the bidding period and exceed a certain minimum amount~\cite{swissgridOkt2015}. 
These requirements reduce the coordination effort for the grid operator by simplifying the scheduling and avoiding a large number of small bidders.
 The minimum bid typically ranges from~1\,MW to 10\,MW, e.g., it is~5\,MW in Switzerland {\cite[Section~2.3.1]{swissgridOkt2015}}, which is too large for a single building~\cite{chapman2016algorithmic}. 
 To still participate in the reserve market, buildings can cooperate~\cite{VrettosTPS2016}, which requires solving a large-scale optimization problem for obtaining a time-constant joint bid. 
 A centralized solution of such a problem faces several issues: $(i)$~a shared investment is necessary to establish and maintain a central computation unit; $(ii)$~solving the resulting large-scale optimization problem can be computationally challenging, which effectively limits the aggregation size; $(iii)$~the building operators need to share sensitive information, such as details on the construction material and occupants' behavior; and $(iv)$~each local configuration change needs to be communicated to the central unit.\looseness=-1

In this paper, we address these issues by using the alternating direction method of multipliers (ADMM)~\cite{glowinski1975approximation,boyd2011distributed} to distribute the original large-scale optimization problem among the individual buildings. 
Each building then solves a localized problem and participates in a global mediation process until a consensus solution is reached. 
The distribution of the computational burden provides natural parallelism and therefore makes it possible to use a simple central processing unit that only handles coordination tasks.
Further, we present a modified communication scheme that avoids a central unit altogether. 
We also show that increased privacy and flexibility follows from the distribution of information, and we demonstrate that we obtain a joint bid satisfying all constraints even if we terminate ADMM prematurely. 
Finally, given that a reserve bid is accepted on the market, we show a reward distribution scheme that reduces provision imbalance within the aggregation, strengthening the group's capability to provide reserves.\looseness=-1
\vspace{-.2\baselineskip}

\subsection*{Related Literature}
Frequency reserve provision by buildings receives considerable attention, both from the power systems community \cite{diekerhof2017hierarchical,Vrettos14ifac,callaway2011achieving, Hao2013, Maasoumy2014, VrettosAndersson_TSE2016,  VrettosTPS2016,chen2013mpc,vrettos2013combined} and the control community \cite{ZhangKamgGoulartLyg_CDC2014, BalandatOldewurtelChenTomlin2014, zhang2015stochastic, GoreckiBitlisliogluStathopoulosJones_2015,   chapman2016algorithmic, ZhangKamgGeorghiouGoulLyg_Aut17, BitlisliogluGoreckiJones_TAC2017,bilgin2016provision,taha2017buildings,kraning2014dynamic}. 
In~\cite{Vrettos14ifac,VrettosTPS2016}, the primary goal is to determine the optimal reserve capacity that an aggregation of buildings can provide while ensuring occupant comfort.
In~\cite{Hao2013, Maasoumy2014}, frequency reserves are provided by modulating the fan speed of air handling units in commercial buildings. 
However, rather than giving a priori guarantees on comfort constraints as in~\cite{Vrettos14ifac,VrettosTPS2016}, the authors in~\cite{Hao2013, Maasoumy2014} test constraint satisfaction a posteriori through experimental trials.
Our work follows~\cite{Vrettos14ifac, VrettosTPS2016}, and we use the mathematical tools for robust constraint satisfaction from~\cite{ZhangKamgGoulartLyg_CDC2014, ZhangKamgGeorghiouGoulLyg_Aut17}. 
In contrast to~\cite{Vrettos14ifac, VrettosTPS2016}, our leading objective is a distributed and decentralized solution of the multi-building problem. 
As opposed to~\cite{vrettos2013combined,weckx2015load}, where voltage regulation affects the reactive power balance, we concentrate on frequency regulation and active power. 
Aside from our focus on market bidding,~\cite{taha2017buildings,kraning2014dynamic,chen2013mpc,bilgin2016provision} analyze the real-time balance in the power grid. In~\cite{taha2017buildings,kraning2014dynamic}, the grid is modeled explicitly, exceeding our work as we purely rely on the reserve demand as requested by the grid operator. Incorporating such grid dynamics can be worthwhile, particularly for spatially separated aggregations. While~\cite{chen2013mpc} guides the demand response based on time-varying energy prices,~\cite{bilgin2016provision} considers a setup where previously agreed reserves are utilized. Hence,~\cite{bilgin2016provision} augments our work particularly well, as we focus on a market where such reserve agreements are made.\looseness=-1

In~\cite{hou2017distributed,verschae2016coordinated,diekerhof2017hierarchical,burger2017generation}, ADMM is studied in the context of demand response.
In~\cite{hou2017distributed}, a distributed procedure coordinates different temperature zones inside a building. In~\cite{verschae2016coordinated}, ADMM exploits the flexibility of several consumers for shaping their combined power intake. In~\cite{diekerhof2017hierarchical}, heterogeneous agents are aggregated. Finally, in~\cite{burger2017generation}, small thermostatically controlled loads are orchestrated to participate in real-time energy markets. 
We differ from these publications by focusing on ahead-of-time reserve bidding while the actual reserve request is still unknown. 
We further set ourselves apart by developing additional tools that improve the conceptual utility of our approach, such as a fully-decentralized computation scheme, the permanent availability of a feasible solution, and a profit allocation that places an incentive to strengthen the group. 
\vspace{-0.2\baselineskip}

\subsection*{Notation}
For \mbox{$x,y\in \mathbb R^n$}, we use $(x;y) = [x^\top\ y^\top]^\top\in \mathbb R^{2n}$. Inequalities between vectors, e.g., $x\leq y$, as well as inclusions, e.g., $x\in [-y,y]$, are understood element-wise. We write $\{x_i\}_i$ for the tuple $(x_1,x_2,\dots)$, where the context determines the range of $i$. We denote the Kronecker product as~$\otimes$. 
We use the identity matrix $I_n$ of dimension $n\times n$.
$0_{n\times m}$ or~$1_{n \times m}$ are $n\times m$ matrices with all elements $0$ or $1$, respectively; we omit the subscripts if the dimension is clear from context. 



\section{Frequency Reserve Provision}\label{sec:setup}
We first introduce the model of a single building and describe the participation in the reserve market. We then show the joint bidding procedure for  aggregations of buildings.

\subsection{Building Model and Reserve Market}
\label{sec:buildingModel}
Similar as in~\cite{OldewurtelenergyBuildings2012,Lehmann_2013,VrettosAndersson_TSE2016}, we consider the building model 
\begin{IEEEeqnarray}{rCl}
\label{eq:buildingmodel}
  x^{k+1} = A x^{k} + B  (u^k + \Delta u^k) + E v^{k},
\end{IEEEeqnarray}
where we use the admissible state $x^k \in \mathcal X \subseteq \mathbb R^n$, the admissible input $(u^k + \Delta u^k) \in \mathcal U \subseteq \mathbb R^m$, and the external disturbance~$v^k\in \mathbb R^q$. We use the discrete time \mbox{$k=1,2,\myDots,N$} with horizon~$N$. 
The state $x^k$ models temperatures in rooms, walls, floors, and ceilings across the building. We use~$\mathcal X$ to~model convex comfort constraints, e.g., to limit room temperatures between $21^\circ C$ and $25^\circ C$.
The system matrix~$A$ describes the temperature diffusion with a thermal resistance-capacitance model as in~\cite[Sec.~4.2]{OldewurtelenergyBuildings2012}.
The actuation (heating, cooling, and ventilation) is described by a nominal input~$u^k$ and a variable term~$\Delta u^k$. The variable term is used later for adjusting the buildings' power consumption. The convex set~$\mathcal U$ models actuation limits and saturation.
Finally, we use~$E v^k$ to describe the effect of external disturbances (occupancy, solar radiation, and ambient temperature). Section~\ref{sec:simulation} shows a simulation study with detailed modeling examples. 
The mechanisms behind building-internal temperature regulation are specifically addressed in~\cite{OldewurtelenergyBuildings2012,Lehmann_2013}. As we focus on the collaboration of different buildings in this paper, we suffice with the abstracted model~\eqref{eq:buildingmodel}.
 Indeed, any linear system~$A,B,E$ with convex constraints~$\mathcal X, \mathcal U$ is compatible with the subsequent analysis.
 We also simplify the setup by assuming perfect prediction for~$v^k$, noting that uncertainty can easily be integrated into our approach by using the robust control methods in\mbox{\cite{Oldewurtel_PhD2011, zhang:schildbach:sturzenegger:morari:13, Darivianakis_CDC2015}}.  With $x \coloneqq (x^2;\myDots;x^{N+1})$, $u \coloneqq (u^1;\myDots;u^{N})$, and similar definitions for~$\Delta u, v$ (ranging from~$1$ to~$N$), we rewrite \eqref{eq:buildingmodel} as \looseness=-1
\begin{IEEEeqnarray}{rCl}
\label{eq:buildingmodelstacked}
x =  \bm A x^1 +  \bm B( u +\Delta   u) +  \bm E v
\end{IEEEeqnarray}
with $x\in\bm{\mathcal X}$ and $( u +\Delta   u) \in \bm{\mathcal U}$. We assemble
$\bm A, \bm B, \bm E, \bm{\mathcal X}, \bm{\mathcal U}$ from $A,  B,  E, \mathcal X, \mathcal U$ as shown in~\cite{Goulart2006}.

\begin{figure}[h!bt]
		\vspace{-0.3\baselineskip}
        \centering
        \includegraphics[width=0.98\columnwidth]{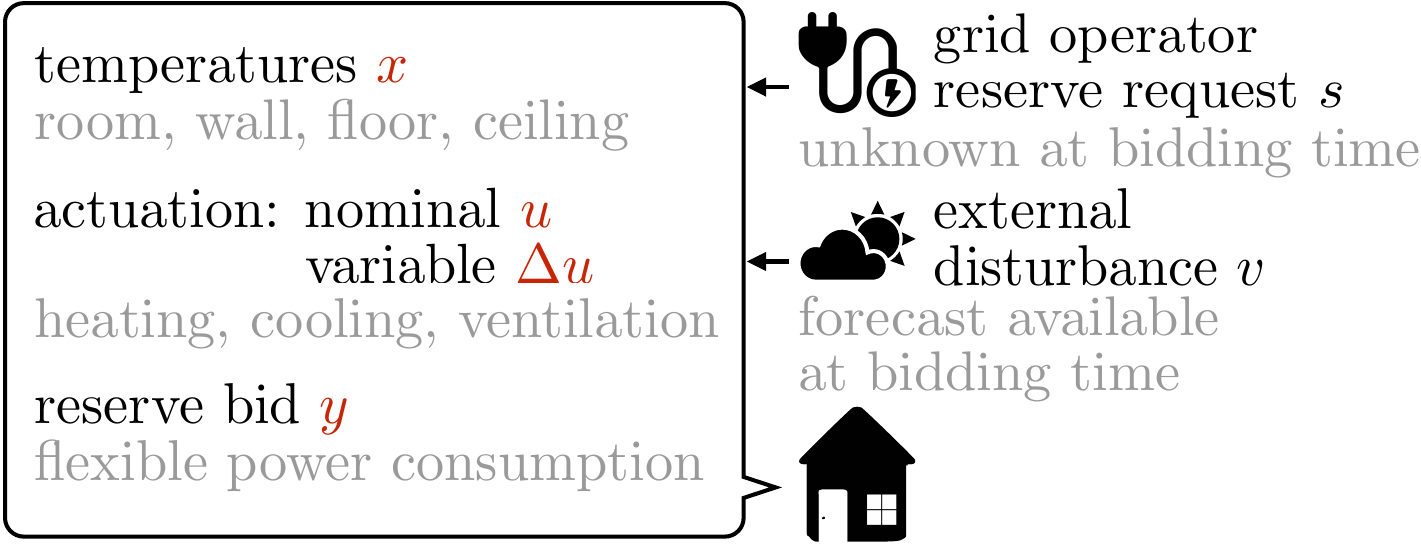}
        \vspace{-0.2cm}
        \caption{Decision variables $x$, $u$, $\Delta u$, $y$, which can be determined by the building, and external variables $s$, $v$, which cannot be influenced.}
        \label{fig:buildingdecision}
\end{figure}

In Switzerland, similar to most of Europe, reserve capacities are procured in a market where reserve providers bid their capacity in weekly or daily blocks. 
We denote the time-varying reserve bid as $y = (y^1;\myDots;y^{N})$, which describes the maximum-possible change in active power consumption, measured in kilowatt. 
We focus on a day-ahead market, and we consider a symmetric reserve range $[-y,y]$, which means that the building needs to shift its consumption equally in both directions. The limitation to symmetric provision follows the current policy in Switzerland~\cite{swissgridOkt2015}; the extension to asymmetric reserves is straightforward. 
We further require a constant bid over the bidding period, i.e., \mbox{$y^1=y^2=\myDots=y^N$}, which is also standard in several countries. The restriction to time-constant reserve simplifies the scheduling for the grid operator and prevents an imbalance over the bidding interval.
Once a bid is accepted, the building receives a reward of the form $p^\top y$, where $p\in\mathbb R^N$ is a price parameter that we assume to be known, e.g., from previous market clearings. 
During actual runtime, the grid operator then requests a reserve demand $s = (s^1;\myDots;s^N)\in [-y,y]$, and the building is obliged to adjust its power consumption accordingly. The challenge in ahead-of-time bidding is choosing~$y$ while~$s$ is still unknown. Fig.~\ref{fig:buildingdecision} summarizes the involved variables.\looseness=-1
\clearpage

\begin{figure}[h!tb]
        \centering
        \includegraphics[width=0.96\columnwidth]{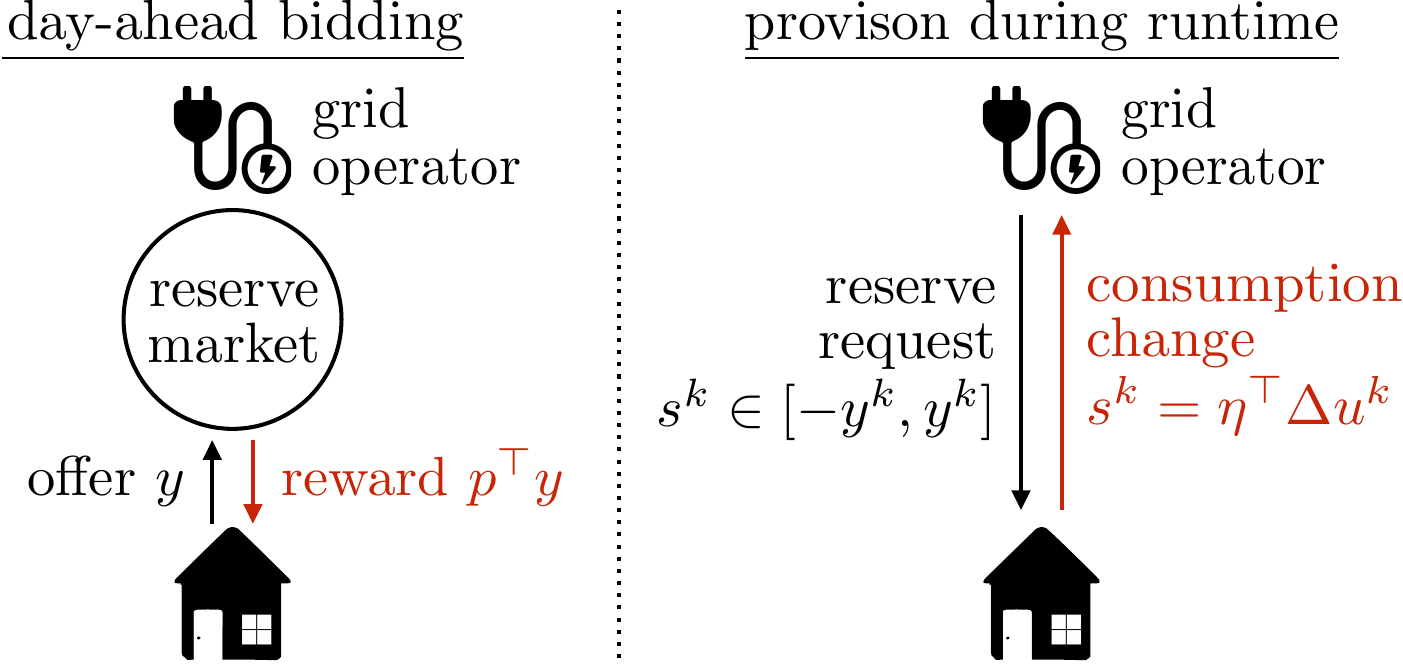}
        \vspace{-0.2cm}
         \caption{Bidding process and reserve provision. 
\emph{Left}:~Ahead of time, the building offers a bid $y$ on the reserve market. If the bid gets accepted, it receives the reward $p^\top y$. The actual reserve request~$s^k$ remains unknown at this point. \emph{Right}: During runtime, the building follows the reserve request $s^k$ within the previously offered range $[-y^k,y^k]$ by adjusting its variable input~$\Delta u^k$.\looseness=-1}
        \label{fig:setup}
\end{figure}

After a bid is accepted at the reserve market, the building follows the request $s$ by appropriately changing its input from~$u$ to~\mbox{$u + \Delta u$}. 
More precisely, following \mbox{\cite{OldewurtelenergyBuildings2012, zhang:schildbach:sturzenegger:morari:13, darivianakis2016power}}, we use the conversion factor $\eta \in\mathbb R^{m}$ for translating input quantities (e.g., heating levels) into power consumption (measured in kilowatt). We choose the variable input~$\Delta u$ such that~$s^k = \eta^\top \Delta u^k$ for all $k$. Fig.~\ref{fig:setup} illustrates the separation between reserve bidding and the actual provision procedure. \looseness=-1

Our approach focuses on the bidding process, particularly on determining a satisfiable and financially profitable bid for the reserve market. Additionally, our bidding method results in a rudimentary procedure to follow the reserve request during runtime, which means that we relate $\Delta u$ to $s$. Beyond the scope of our paper, our provision procedure can be complemented with more sophisticated and iterative reserve tracking techniques as in~\cite{taha2017buildings,kraning2014dynamic,bilgin2016provision}, which can exploit newly-available information or incorporate grid models to fine-tune the supply.\looseness=-1

\subsection{Individual Bidding Problem}
Typically, buildings need to increase their nominal consumption $u$ for being able to react to symmetric reserve requests~\cite{VrettosTPS2016}.
The building's objective is to optimally balance the provision reward $p^\top  y$ and the nominal electricity cost $c^\top  u$, i.e., it minimizes its overall cost $c^\top u - p^\top y$. Here, the vector $c$ is composed of subvectors $c_k = \tilde c_k \eta$, where~$\tilde c_k\in \mathbb R$ is the electricity price, and $\eta$, as introduced before, converts input signals to power consumption. 
Since the reserve request~$s\in\mathbb R^N$ is unknown at the time of bidding, 
we introduce recourse into the formulation by choosing $u$ and $\Delta u$ to be closed-loop policies 
\begin{IEEEeqnarray}{rClrCl}
\IEEEyesnumber\eqlabel{eq:causalfunDef} \IEEEyessubnumber*
	 u(\cdot)&:&\mathbb R^N \rightarrow \mathbb R^{N\hspace{-0.2ex}m} \quad & s &\mapsto&  u\\
	 \Delta u(\cdot)&:&\mathbb R^N \rightarrow \mathbb R^{N\hspace{-0.2ex}m} \quad & s &\mapsto&  \Delta u,
\end{IEEEeqnarray}
which are parametrized functions that determine the actual input only when~$s$ becomes known during runtime \cite{Scokaert1998,LofbergPhD2003, Goulart2006}. In other words, instead of determining the static values~$u, \Delta u$ that prepare for the worst-case request~$s$, we try to find the optimal functional relation of~$u, \Delta u$ and~$s$. Compared to the conventional strategy, closed-loop policies offer a larger degree of freedom and therefore  potentially less conservative reserve bids~\cite{zhang2015stochastic}. We choose the policies to be causal such that~$u^k, \Delta u^k$ only depend on~$s^1,\dots,s^k$. The building-individual reserve provision problem is\looseness=-1
\begin{IEEEeqnarray}{rCl}
  \IEEEyesnumber\eqlabel{eq:individualProgram} 
  \min_{ u(\cdot), y} &\quad&  c^\top  u(0) -  p^\top  y \IEEEyessubnumber \label{individualObjective}\\[-0.2ex]
  \text{s.t.} && \left( u(\cdot), y\right) \in \mathcal C \IEEEyessubnumber\label{eq:consumerConstraint}\\
  && y^1 = y^2 =\cdots=y^N.\IEEEyessubnumber\label{eq:yconst}\vspace{-0.5ex}
\end{IEEEeqnarray}
 
For simplicity, we use the nominal electricity cost~$c^\top  u(0)$, which is based on the case $s=0$ where the grid operator does not request any reserves. Due to~$\Delta u(0)=0$, the variable input is not part of the objective.
Alternatively, it is possible to use the expected cost~$\mathbb E\left[ c^\top \left( u( s) + \Delta u( s) \right) \right]$, provided we have a probabilistic description of the demanded reserve~$s$~\cite{zhang2015stochastic}. 
Further,
we choose not to model the minimum reserve constraint~\mbox{$y^k\geq y_\text{min}$}. Having such a constraint may lead to overall increasing costs compared to not providing reserves at all. 
Instead, we purely maximize financial profit, and we assume that the resulting bid is only placed at the market if it exceeds the minimum size. 
Omitting the minimum reserve constraint also benefits the construction of the distributed setup, as we will clarify in Remark~\ref{rem:whatifminbid}.  
 The set~$\mathcal C$ in \eqref{eq:consumerConstraint} captures the building dynamics and the reserve provisioning mechanism as introduced before. It encapsulates building-internal information, e.g., the state $x$, which in consequence of~\eqref{eq:causalfunDef} also becomes a policy.
 We define $\mathcal C$ as\looseness=-1
 \begin{IEEEeqnarray}{l}
  \IEEEyesnumber \label{eq:consumerConstraintExample} 
\mathcal C = \big\{\left( u(\cdot),  y \right) \in \left(\mathcal F \times \mathbb R^N \right)  \big|\, \exists \left(\Delta  u(\cdot),  x(\cdot)\right)\in \mathcal F \times \mathcal F \text{ s.t.}\nonumber \\
  \left. \begin{array}{r}
   x( s) =  \bm A x^1 +  \bm B[ u(s) +{\Delta   u( s)}] +  \bm E v\\
  \left( x( s), u( s) + \Delta  u( s) \right) \in \bm{\mathcal X} \times \bm{\mathcal U}\\
   s =  (I_N \otimes \eta^\top) \Delta u( s)
\end{array} \hspace{-1mm} \right\} \forall  s \in [-  y,  y]\, \big\},\nonumber\\[-2.5ex] \vspace{-0.5ex}
\end{IEEEeqnarray}
where $\mathcal F$, in the most general case, is the infinite-dimensional space of causal functions of type~\eqref{eq:causalfunDef}.
Due to the optimization over policies and the infinite number of constraints, solving~\eqref{eq:individualProgram} is intractable in general. However, if we restrict~$\mathcal F$ to contain linear functions only, we can reformulate~\eqref{eq:individualProgram} as a convex optimization problem with a finite number of constraints; see 
Appendix~\ref{app:defCi} and~\cite[Section 4]{ZhangKamgGeorghiouGoulLyg_Aut17} for an example. Given such an appropriate choice of~$\mathcal F$, the constraint set~$\mathcal C$ and therefore the individual bidding problem~\eqref{eq:individualProgram} is finite-dimensional and convex. Indeed, finite dimensionality and convexity of $\mathcal C$ are the only requirements for the analysis below. Hence, our approach is compatible with different policy spaces and building models. Further, it can be used for plug-in electric vehicles, house-mounted batteries, and any energy storage device that has linear dynamics $\bm A, \bm B, \bm E$ and convex constraints $\bm{\mathcal X}, \bm{\mathcal U}$.\looseness=-1

\begin{remark} 
In some markets, the price~$p$ is also a decision variable, which increases the computational cost significantly. To keep the computation simple, we can solve~\eqref{eq:individualProgram} for different values of~$p$, which results in a price-reserve curve that is bid at the market. This strategy is known as conditional bidding.\looseness=-1
\end{remark}

\subsection{Aggregated Bidding Problem}
\label{sec:aggregatedBiddingProblem}
On their own, buildings are typically not able to meet the required minimum bid size~\cite{VrettosAndersson_TSE2016}. To overcome this, buildings can collaborate in an aggregation.
We consider $M$ buildings, indexed by~\mbox{$b=1,\myDots,M$}, and we denote the reserve bids as~$y_b$, where each building's bid has the form \mbox{$y_b=(y_{b}^1;\myDots; y_{b}^N)$}. Similarly, we use~$u_b$, $c_b$, and~$\mathcal C_b$ for the individual  input, cost function, and constraint set. Buildings are compatible as long as they agree on the bidding horizon $N$, which means that they can use different models in~$\mathcal C_b$. We use non-calligraphic uppercase symbols for quantities associated with the entire aggregation, e.g., we define the aggregated reserve \mbox{$Y\coloneqq \sum\nolimits_{b}  y_b$}, where \mbox{$Y = (Y^1;\myDots;Y^N)$}. The aggregated bidding problem is\looseness=-1
\begin{IEEEeqnarray}{cCl}
 \IEEEyesnumber\eqlabel{eq:aggregatedProgram} 
  \min_{ Y,\{u_b(\cdot), y_b\}_{b}} &\quad& \sum\nolimits_{b=1}^M  c_b^\top  u_b(0) -  p^\top Y \IEEEyessubnumber\\
  \text{\,\,\,s.t.} && \left( u_b(\cdot), y_b \right) \in \mathcal C_b \quad \forall b\IEEEyessubnumber\label{eq:inC}\\
  && Y = \sum\nolimits_{b=1}^M  y_b \IEEEyessubnumber\label{eq:Ydef}\\
  && Y^1 = Y^2 = \cdots = Y^N. \IEEEyessubnumber\label{eq:Yconst}
\end{IEEEeqnarray}

In the aggregated problem, the time-constant constraint~\eqref{eq:Yconst} applies to the aggregated bid~$Y$, which gives each building the freedom of time-varying bids~$y_b$ as long as their sum remains constant. This additional freedom increases the joint bid, as compared to the trivial strategy of pooling individual bids obtained with~\eqref{eq:individualProgram}. The potential increase in~$Y$ comes at the cost of requiring a deeper collaboration between the buildings, as each building has to shape its offer~$y_b$ while considering the bids of all other buildings. If the time-constant constraints~\eqref{eq:yconst} and~\eqref{eq:Yconst} would not be present, then~\eqref{eq:aggregatedProgram} would simply decompose into separate instances of~\eqref{eq:individualProgram}, i.e., there would be no need for collaboration.  When the buildings agree on a joint bid $Y$ that is also accepted at the reserve market, then the grid operator can request his demand~$s$ from the aggregation during runtime. Each building then changes its consumption proportionally to its share~$y_b$ of the original offer~$Y$.
 Beyond the scope of our paper, it is possible to successively reallocate the provision responsibilities within the aggregation during runtime, which can be especially useful when taking a grid model and therefore local imbalances into account~\cite{taha2017buildings,kraning2014dynamic}.\looseness=-1

A direct solution of~\eqref{eq:aggregatedProgram} has several downsides:
$(i)$~Due to the large dimensionality, we require a powerful central processing unit, which demands financial resources and restricts the aggregation size through limited computational capabilities. 
$(ii)$~The buildings need to disclose $\mathcal C_b$, potentially including sensitive information, which also implies that they cannot change $\mathcal C_b$ without informing the central unit.
$(iii)$~The buildings delegate the decision process, i.e., they have to rely on the central unit for finding a good solution to~\eqref{eq:aggregatedProgram}.\looseness=-1

\section{Aggregated Bidding with ADMM}
\label{sec:admm}
We use ADMM to decompose~\eqref{eq:aggregatedProgram} into a sequence of smaller problems, which results in conceptual benefits over a direct solution approach. We rely on ADMM due to its mild convergence requirements, e.g., we do not need strong convexity.\looseness=-1

\subsection{General ADMM Formulation}
ADMM is a first-order optimization method that solves problems where the decision variables are partitioned into two groups. Typically, these problems have the form 
\begin{IEEEeqnarray}{rCl}
  \label{eq:standardADMM}
  \min\nolimits_{z, w}  \left\{\, f(z) + g(w) \,\text{~s.t.~} \,\mathcal A(z) + \mathcal B(w) = 0\, \right\}
\end{IEEEeqnarray} 
with real-valued functions $f, g$ and affine operators $\mathcal A, \mathcal B$.\looseness=-1

As shown in Alg.~\ref{alg:admm}, ADMM solves~\eqref{eq:standardADMM} by alternatingly searching for a saddle point of the augmented Lagrangian
\begin{IEEEeqnarray}{rCl}
\label{eq:Langrangiandef}
  \mathcal L_\rho \left( z, w, \lambda \right) &\coloneqq& 
  f(z)
  + g(w)
  + { \lambda^\top \epsilon} 
  + \tfrac{\rho}{2} \left\| \epsilon  \right\|^2_2,
\end{IEEEeqnarray}
where 
$\lambda$ is the Lagrange multiplier, $\epsilon \coloneqq \mathcal A(z) + \mathcal B(w)$ is a residual, and $\rho>0$ is a user-chosen penalty parameter. 
For any initialization of $(z, w,\lambda)$, and any positive value of~$\rho$, ADMM converges to a solution of~\eqref{eq:standardADMM}, given that a solution exists, and~$f,g$ have a nonempty, closed, and convex epigraph~\cite{boyd2011distributed}.

\begin{algorithm}[h]
\caption{ADMM}\label{alg:admm}
\begin{algorithmic}[1]
\vspace{0.2em}

\Statex \hspace{-1.05em}\textbf{repeat}
\vspace{0.1em}
\,
\begin{minipage}[t]{\linewidth}
\begin{algorithmic}[1]
\State $z\, \leftarrow \arg \min\limits_{z}{\, \mathcal L_\rho \left( z, w, \lambda \right) }$
\vspace{0.2em}
\State $w \leftarrow \arg \min\limits_{w}{ \mathcal L_\rho \left( z, w, \lambda \right) }$
\vspace{0.2em}
\State $\lambda \,\leftarrow \lambda + \rho  {\tfrac{\partial}{\partial \lambda} \mathcal L_\rho \left( z,w, \lambda \right) }$
\end{algorithmic}
\end{minipage}
\vspace{0.01em}
\Statex \hspace{-1em}\textbf{until} {satisfaction of a stopping criterion}
\end{algorithmic}
\end{algorithm}
\vspace{-1\baselineskip}

\subsection{Aggregated Bidding with ADMM}
\label{sec:aggregatedBiddingviaADMM}
For writing~\eqref{eq:aggregatedProgram} in the form of~\eqref{eq:standardADMM}, we first need to relate the decision variables~$Y,\{u_b(\cdot), y_b\}_{b}$ to the two-component splitting~$z,w$. Further, we need to associate each constraint in~\eqref{eq:aggregatedProgram} to one of the two groups, and we need to decide on a coupling relation $\mathcal A$, $\mathcal B$. Our goal is to associate~$w$ with the aggregation-specific variable $Y$, such that~$z$ only contains building-individual variables $\{u_b(\cdot), y_b\}_{b}$, making the $z$-update in Alg.~\ref{alg:admm} decompose among buildings. When pursuing such a splitting, we struggle with assigning~\eqref{eq:Ydef} to one side, as it connects the aggregation-specific and building-individual variables. 
We overcome this difficulty by making~\eqref{eq:Ydef} depend on the auxiliary decision variable $\{\bar y_{b}\}_b$, which we constrain to be equal to~$\{y_b\}_b$. 
Later on, it will become clear that having these extra variable does not increase the computational burden.
With $\{\bar y_{b}\}_b$, we then can separate between the aggregation-specific variables \mbox{${w} \coloneqq \big({Y},\{\bar y_{b}\}_b\big)$}, and the building-individual variables \mbox{${z} \coloneqq \{u_{b}(\cdot),y_{b}\}_b$}, where we use $\{\bar y_{b}\}_b=\{y_b\}_b$ as coupling relation. The resulting problem partition, shown in~Fig.~\ref{fig:splitting}, is an equivalent reformulation of~\eqref{eq:aggregatedProgram}.
\definecolor{lightblue}{RGB}{3,101,192}
\definecolor{darkred}{RGB}{200,37,5}
\begin{figure}[hbt]
\vspace{-0.5\baselineskip}
      \centering
      \begin{tikzpicture}
      \coordinate(middleBottom) at (0,-.5);



      \coordinate(problemFix) at ($(middleBottom)+(-.6,-.1)$);
      \node[align=center, anchor = east] at ($(problemFix)+(-1,-.8)$)
      	{$\min\limits_{ \substack{\color{lightblue} \{u_b(\cdot), y_b\}_{b}\color{black}, \\ \color{darkred} Y,\{\bar y_b\}_{b}} }$};
      \node[align=center, anchor = east] at ($(problemFix)+(-1.46,-1.4)$)
      	{s.t.};

      \node[align=left, anchor = west] at ($(problemFix)+(-.8,-.8)$)
      {$\sum\nolimits_{b}  c_b^\top  {\color{lightblue}u_b(0)} \hspace{1.1cm} -  p^\top {\color{darkred}Y}$};
      \node[align=left, anchor = west] at ($(problemFix)+(-.8,-1.4)$)
      {$( {\color{lightblue}u_b(\cdot)}, {\color{lightblue}y_b}) \in \mathcal C_b \,\, \forall b$};
      \node[align=left, anchor = west] at ($(problemFix)+(-.8,-2)$)
      {${\color{darkred}Y} = \sum\nolimits_{b}  {\color{darkred}\bar y_b}$};
      \node[align=left, anchor = west] at ($(problemFix)+(-.8,-2.5)$)
      {${\color{darkred}Y^1} = {\color{darkred}Y^2} = \cdots = {\color{darkred}Y^N}$};
      \node[align=left, anchor = west] at ($(problemFix)+(-.8,-3)$)
      {${\color{darkred}\bar y_{b}} - {\color{lightblue} y_{b}} =0 \,\,\, \forall b$};
      
      \coordinate(indTL) at ($(problemFix)+(-.8,-.5)$);
      \coordinate(indTR) at ($(problemFix)+(2.1,-.5)$);
      \coordinate(indBL) at ($(problemFix)+(-.8,-1.65)$);
      \coordinate(indBR) at ($(problemFix)+(2.1,-1.65)$);
      \draw[-, solid, color = lightblue] (indTL)-- (indTR);
      \draw[-, solid, color = lightblue] (indTR)-- (indBR);
      \draw[-, solid, color = lightblue] (indBR)-- (indBL);
      \draw[-, solid, color = lightblue] (indTL)-- (indBL);

      \coordinate(agrTL) at ($(problemFix)+(-.8,-1.75)$);
      \coordinate(agrTR) at ($(problemFix)+(3.4,-.5)$);
      \coordinate(agrBL) at ($(problemFix)+(-.8,-2.75)$);
      \coordinate(agrBR) at ($(problemFix)+(3.4,-2.75)$);
      \coordinate(agrMB) at ($(problemFix)+(2.2,-1.75)$);
      \coordinate(agrMT) at ($(problemFix)+(2.2,-.5)$);
      \draw[-, solid, color = darkred] (agrTL)-- (agrMB);
      \draw[-, solid, color = darkred] (agrMB)-- (agrMT);
      \draw[-, solid, color = darkred] (agrMT)-- (agrTR);
      \draw[-, solid, color = darkred] (agrTR)-- (agrBR);
      \draw[-, solid, color = darkred] (agrBR)-- (agrBL);
      \draw[-, solid, color = darkred] (agrTL)-- (agrBL);
      \coordinate(aggrRelCoord) at ($(problemFix)+(3.7,-1.75)$);
      \node[align=left, anchor = west] at (aggrRelCoord) {\color{darkred} aggregation-\\ \color{darkred}specific};

      \coordinate(builRelCoord) at ($(problemFix)+(3.7,-.8)$);
      \node[align=left, anchor = west] at (builRelCoord) {\color{lightblue} building-\\ \color{lightblue}individual};

\end{tikzpicture}
\vspace{-.5\baselineskip}
\caption{Separation of~\eqref{eq:aggregatedProgram} into building-individual and aggregation-specific terms, where we use the auxiliary variable $\bar y_{b} =  y_{b}$ for $b=1,\myDots, M$.}\label{fig:splitting}
\end{figure}

By comparing the partitioned problem in Fig.~\ref{fig:splitting} to the general ADMM problem~\eqref{eq:standardADMM}, we obtain the definitions
\begin{IEEEeqnarray}{rCl}
\IEEEyesnumber \label{eq:splittingDefinition} \IEEEyessubnumber*
  f(z) &\coloneqq&
   \sum\nolimits_{b} \big(c_b^\top u_b(0) + \mathcal I_{\mathcal C_b}\left(u_b(\cdot),y_b\right) \hspace{-2pt}\big) \\
  g(w) 
  &\coloneqq& \mathcal I_{Y = \sum_{b} \bar y_b}(Y, \{\bar y_{b}\}_b) + \mathcal I_{Y^1=\cdots = Y^N} (Y) - p^\top Y \quad \quad \\
  \mathcal A(z) \hspace{-.16667em} &+&\hspace{-.16667em}\mathcal B(w) = 0 \Leftrightarrow \{\bar y_{b}\}_b =  \{ y_{b}\}_b,
\end{IEEEeqnarray}
where indicator functions are denoted as $\mathcal I(\cdot)$. Inserted into~\eqref{eq:Langrangiandef} and Alg.~\ref{alg:admm}, these definitions lead to our proposed ADMM formulation, which we show in Alg.~\ref{alg:ouradmm}.

\newpage
\phantom{v}
\vspace{-1.4\baselineskip}

\begin{algorithm}[h]
\caption{ADMM for Reserve Provision}\label{alg:ouradmm}
\begin{algorithmic}[1]
\Statex \hspace{-1em}\textbf{repeat}

\algrenewcommand{\alglinenumber}[1]{\hspace{1.4ex}\footnotesize\circled{1}\hspace{-0.4ex}}
\State building-individual reserve proposition $\left(b=1,\myDots,M\right)$
\vspace{0.5mm}
\Statex $\begin{array}{rl}
(u_b(\cdot),y_b) \leftarrow \argmin \limits_{u_b(\cdot),y_b}
&\hspace{-0.6em} c_b^\top u_b(0) - \lambda_b^\top y_b + \tfrac{\rho}{2}\| \bar  y_b - y_b \|_2^2 \\
\text{s.t.}   &\hspace{-0.6em} \left(u_b(\cdot),y_b \right) \in \mathcal C_b
\end{array}$ \vspace{-0.3em}

\algrenewcommand{\alglinenumber}[1]{\hspace{1.4ex}\footnotesize\circled{2}\hspace{-0.4ex}}
\State aggregation step 
\vspace{-1.5mm}
\Statex \hspace{-3mm}$\hspace{1.6mm} \begin{array}{rl} (Y, \{\bar y_{b}\}_b) \leftarrow \argmin\limits_{Y, \{\bar y_b\}_b}
&\hspace{-0.6em} \sum \limits_{b=1}^M \left( \lambda_b^\top \bar y_b + \tfrac{\rho}{2}\| \bar y_b -  y_b \|_2^2 \right) - p^\top Y\\
\text{s.t.}   	&\hspace{-0.6em} Y = \sum_{b} \bar y_b \\[.4ex]
        		&\hspace{-0.6em} Y^1 = Y^2 = \cdots = Y^N 
\end{array}$ \vspace{-.5mm}

\algrenewcommand{\alglinenumber}[1]{\hspace{1.4ex}\footnotesize\circled{3}\hspace{-0.4ex}}
\State Lagrangian mediation $\left(b=1,\myDots,M\right)$ 
\Statex \hspace{1.6mm} $\lambda_b \leftarrow \lambda_b + \rho(\bar y_b - y_b)$\vspace{0.2em}

\Statex \hspace{-1em}\textbf{until} {satisfaction of the stopping criterion}
\end{algorithmic}
\end{algorithm}

 Alg.~\ref{alg:ouradmm} alternates between a building-individual step~\circled{1}, an aggregation-specific step~\circled{2}, and Lagrangian mediation~\circled{3}. 
For the initialization, we can set $Y$, $\bar y_b$, and $\lambda_b$ to zero. Prop.~\ref{prop:ouradmmconvergence}, proven in Appendix~\ref{app:proofouradmmconvergence}, constitutes the use of Alg.~\ref{alg:ouradmm}.\looseness=-1

\begin{proposition}
\label{prop:ouradmmconvergence}
If $\mathcal C_b$ is nonempty, finite-dimensional, and convex for each building $b$, then Alg.~\ref{alg:ouradmm} converges to an optimal solution of the aggregated bidding problem~\eqref{eq:aggregatedProgram}.
\end{proposition}
 
\begin{figure}[h!bt]
        \centering
        \includegraphics[width=0.70\columnwidth]{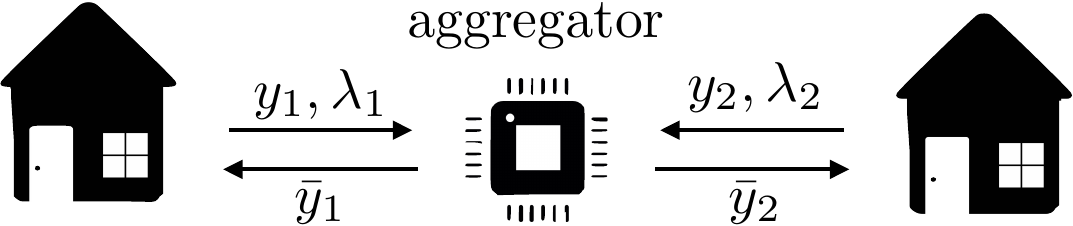}
        \vspace{-0.2cm}
        \caption{Information exchange between the buildings and the aggregator within one iteration of Alg.~\ref{alg:ouradmm}.  Illustrated is the case of $M=2$.}
        \label{fig:communication}
\end{figure}

We interpret Alg.~\ref{alg:ouradmm} as a negotiation process between the buildings and the aggregator.  In~\circled{1}, each building proposes a reserve capacity~$y_b$ that maximizes its profit.
The aggregator, which receives these propositions, responds with the desired reserve amount~$\bar y_b$, chosen in~\circled{2} for achieving a time-constant total reserve~$Y$. The Lagrangian update~\circled{3} forms a price mechanism that mediates between proposition~$y_b$ and demand~$\bar y_b$ to reach a consensus. 
Fig.~\ref{fig:communication} illustrates the information exchange.\looseness=-1

The ADMM-based approach shows conceptual advantages when compared to a centralized solution of~\eqref{eq:aggregatedProgram}. Due to the building-individual execution of~\circled{1}, the constraint sets~$\mathcal C_b$ remain private, and the buildings can change them self-organized.  Also, the execution  parallelizes naturally. Other benefits of Alg.~\ref{alg:ouradmm} derive from the simplicity of its aggregation step, which is described in the following proposition.
\begin{proposition}
\label{prop:analyticaggregatorupdate} 
The aggregation step~\circled{2} has the solution\looseness=-1
\begin{IEEEeqnarray}{rCl}
\IEEEyesnumber \label{eq:analyticaggregatorupdate}  \IEEEyessubnumber*
\Omega &\coloneqq& \tfrac{1}{M}\sum\nolimits_{b=1}^M (\rho y_b - \lambda_b)\label{eq:Omegaresult}\\
  Y &=& \tfrac{M}{\rho N} 1_{N\times N} \left( \Omega + p \right)  \label{eq:Yresult}
  \\
  \bar y_b &=& \tfrac{1}{\rho}\left(\rho y_b  -\lambda_b - \Omega\right) 
+ \tfrac{1}{M}Y     
   \quad \forall b,  \label{eq:ybresult}
\end{IEEEeqnarray}
where $\Omega$ is called the aggregate.
\end{proposition}%
\iflongversion
 Prop.~\ref{prop:analyticaggregatorupdate}, which we prove in Appendix~\ref{app:proofPropAnalytic}, makes it possible to replace the aggregation step with a set of closed-form algebraic equations. 
This signifies that the computational complexity of~\circled{2} is negligible compared to the first (parallel) algorithm step. As a consequence, the increased computational burden for large aggregations is fully absorbed by the equally growing level of parallelism. The ADMM-based framework also lowers the barriers to integrate new buildings into an existing aggregation. The reason is the similarity between the individual problem~\eqref{eq:individualProgram} and the building-individual step \circled{1} in Alg.~\ref{alg:ouradmm}, i.e., new buildings that switch from individual bidding practically keep performing the same optimization as before.\looseness=-1


\begin{remark}
\label{rem:whatifminbid}
	If we model a minimum bid constraint $y^k\geq y_\text{min}$ in~\eqref{eq:individualProgram}, and accordingly $Y^k\geq y_\text{min}$ in~\eqref{eq:aggregatedProgram}, then the aggregation step contains inequality constraints and therefore does not have a closed-form solution anymore.
\end{remark}

\subsection{Decentralized Computation}
Alg.~\ref{alg:ouradmm} is flexible against configuration changes, avoids concentration of private information, and behaves well from a computational perspective. Still, the algorithm requires the acquisition and maintenance of a central computation facility, which poses an organizational burden. To overcome this, we recognize that in the aggregation step~\eqref{eq:analyticaggregatorupdate}, only the aggregate~$\Omega$ requires aggregation-global information. More specifically, we notice that when~$\Omega$ is available, then each building can evaluate the previously centralized operations \eqref{eq:Yresult} and~\eqref{eq:ybresult} individually. Towards utilizing this, we define the partial aggregate
\begin{IEEEeqnarray}{rCl}
\Omega_i \coloneqq \tfrac{1}{M}\sum\limits_{b=1}^i (\rho y_b - \lambda_b),
\end{IEEEeqnarray}
which we can accumulate by serially handing it from building to building. Each building $b$ then adds its own contribution~$\tfrac{1}{M}(\rho y_b - \lambda_b)$, assuming that the aggregation size $M$ is known.  Once we have reached $\Omega = \Omega_M$, the buildings circulate the result back for making it known to their neighbors. Fig.~\ref{fig:comm_decentral} illustrates the communication procedure for the case of four buildings. In general, we can use any connected communication graph that covers the aggregation.
\begin{figure}[h!bt]
        \centering
        \includegraphics[width=0.75\columnwidth]{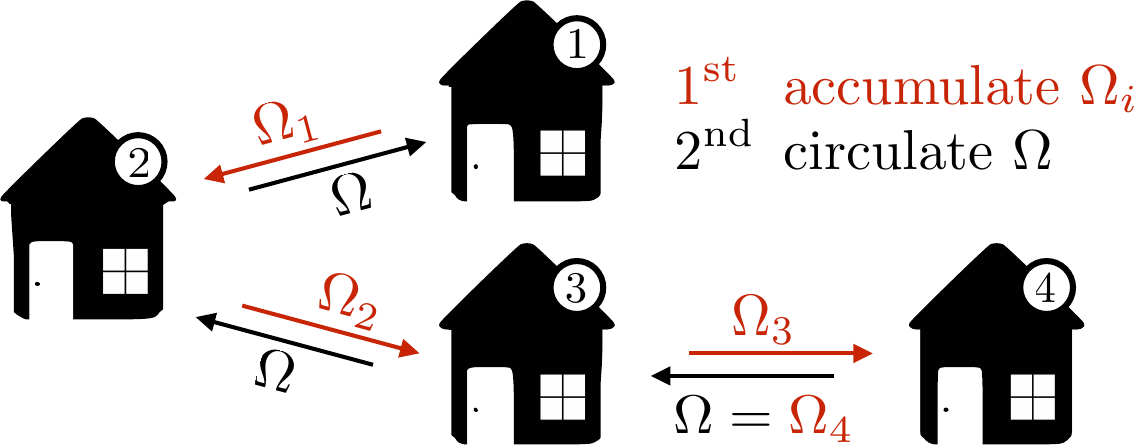}
        \vspace{-0.2cm}
        \caption{Decentralized information exchange to make the aggregate $\Omega$ available to all buildings.  Illustrated is the case of $M=4$.}
        \label{fig:comm_decentral}
\end{figure}

Alg.~\ref{alg:ouradmmDecentral} summarizes the resulting decentralized formulation. The algorithm is obtained by applying the previously described decentralization procedure to Alg.~\ref{alg:ouradmm}, and by rotating its steps from \circled{1}-\circled{2}-\circled{3}  to \circled{2}-\circled{3}-\circled{1}. The result is a grouped building-individual step~\circled{i} and a serial communication step \circled{ii}. Alg.~\ref{alg:ouradmmDecentral} is an equivalent reformulation of Alg.~\ref{alg:ouradmm}, i.e., both algorithms produce the same sequence of iterates when initialized the same. Therefore, all properties that are discussed in this paper equally apply to both algorithms. In subsequent sections, we again focus on Alg.~\ref{alg:ouradmm}, utilizing its more explicit aggregator-building composition to simplify the arguments that follow.\looseness=-1

\newpage
\phantom{v}
\vspace{-1.5\baselineskip}
\begin{figure}[htb]
\vspace{-\baselineskip}
\begin{algorithm}[H]
\caption{Decentralized ADMM for Reserve Provision}\label{alg:ouradmmDecentral}
\begin{algorithmic}[1]
\Statex \hspace{-1em}\textbf{repeat}

\algrenewcommand{\alglinenumber}[1]{\hspace{1.4ex}\footnotesize\circled{i}\hspace{-0.4ex}}
\State building-individual updates $\left(b=1,\myDots,M\right)$\quad 
\vspace{.3em}
\Statex \hspace{-.2em}$\begin{array}{rcll}
\bar y_b &\hspace{-.6em}\leftarrow&\hspace{-.6em} 
\tfrac{1}{\rho} \rlap{$\left(\rho y_b  -\lambda_b - \Omega\right) 
+ \tfrac{1}{\rho N} 1_{N\times N} \left( \Omega + p \right)$}\\[1.4ex]
\lambda_b &\hspace{-.6em}\leftarrow&\hspace{-.6em} 
\lambda_b + \rlap{$\rho(\bar y_b - y_b)$}\\[1.1ex]
(u_b(\cdot),y_b) &\hspace{-.6em}\leftarrow&\hspace{-.6em} 
\argmin \limits_{u_b(\cdot), y_b}&\hspace{-0.4em} c_b^\top u_b(0) - \lambda_b^\top y_b + \tfrac{\rho}{2}\| \bar  y_b - y_b \|_2^2 \\
 && \hfill \text{s.t.}   &\hspace{-0.4em} \left(u_b(\cdot),y_b \right) \in \mathcal C_b
\end{array}$  \vspace{0.1em}

\algrenewcommand{\alglinenumber}[1]{\hspace{1.4ex}\footnotesize\circled{ii}\hspace{-0.4ex}}
\State serial communication 
\Statex \textit{accumulate} $\Omega_i$ \textit{and circulate} $\Omega$ \textit{as shown in Fig.~\ref{fig:comm_decentral}} \vspace{0.5em}

\Statex \hspace{-1em}\textbf{until} {satisfaction of the stopping criterion}
\Statex \hspace{-1em}\textbf{end} obtain the aggregated bid $Y$ with~\eqref{eq:Yresult}
\end{algorithmic}
\end{algorithm}
\vspace{-2\baselineskip}
\end{figure}

\subsection{Feasible Extraction}
\label{sec:feasExtract}
In Alg.~\ref{alg:ouradmm}, the partial results of~\circled{1} and~\circled{2} satisfy their respective constraints; however are not compatible as it is \mbox{$y_b \neq \bar y_b$} before convergence. In other words, ADMM only finds a solution that jointly satisfies all constraints in the limit, potentially requiring a large number of iterations~\cite{rey2016}. Here, we overcome this difficulty with a method for extracting a jointly feasible solution after any number of ADMM iterations. This extraction procedure has great value in practice as it makes it possible to terminate the algorithm anytime, e.g., when we reach a given time limit. The following proposition describes the procedure. \looseness=-1 
\begin{proposition}
\label{prop:feasextract}
Given $\left\{ u_b(\cdot),y_b \right\}_{ b}$ from any iteration of Alg.~\ref{alg:ouradmm}, we obtain a suboptimal solution $\left(Y^F,\left\{ u_b^F(\cdot),y_b^F\right\}_{ b}\right)$ that satisfies all constraints in~\eqref{eq:aggregatedProgram} with 
\begin{IEEEeqnarray}{rCl}
\IEEEyesnumber \IEEEyessubnumber*
  Y^F &\coloneqq& 1_{N\times 1}\min_k \Big(\sum\nolimits_{b=1}^M y_{b}^k \Big) \label{eq:createfeasY}\\
  {\left(y^k_{b}\right)}^F &\coloneqq& \frac {{\left(Y^k\right)}^F}{\sum\nolimits_{j=1}^M y_{j}^k} \, y_{b}^k\quad \forall\,b,k \label{eq:createfeasy}\\
  u^F_{b}(\cdot) &\coloneqq& \argmin_{u_{b}(\cdot)} \left\{ c_b^\top u_{b}(0) \,
	  \text{\,s.t.\,} \left(u_{b}(\cdot),y^F_{b}\right) \in \mathcal C_b \right\} \,\, \forall b.\quad\quad\,\, \label{eq:createfeasu}
\end{IEEEeqnarray}
\end{proposition}

 We show the proof of Prop.~\ref{prop:feasextract} in Appendix~\ref{app:proofPropFeasible}. 
According to~\eqref{eq:createfeasY}, we obtain the feasible combined bid~$Y^F$ by clipping the accumulated reserve offer to its lowest value over time. With~\eqref{eq:createfeasy}, the individual reserve contributions~$y_b$ are scaled down to~$y_b^F$, which then satisfy~$Y^F = \sum\nolimits_{b} y^F_{b}$. In~\eqref{eq:createfeasu}, we reoptimize the input. Due to~$y_b^F\leq y_b$ the previous policy~$u_b(\cdot)$ also remains feasible; however is less efficient.
The more algorithm iterations we perform before extracting a solution, the lower the overall cost $\sum_b c_b^\top u_b^F(0) - p^\top Y^F$. In Section~\ref{sec:stopppingCriteria}, we empirically analyze how many iterations are necessary to come reasonably close to the optimum. 

\subsection{Reward Distribution}
When the grid operator accepts a reserve bid $Y$, the aggregation receives the reward $R \coloneqq p^\top Y$. We can distribute this reward in proportion to the offered reserves, i.e., building $b$ gets\looseness=-1
\begin{IEEEeqnarray}{rCl}
  \IEEEyesnumber \label{eq:reservereward} 
r_{b} \coloneqq p^\top y_b. 
\end{IEEEeqnarray}
It is $R = \sum_b r_b$, signifying that we do not distribute more or less than the total reward.
A drawback of this approach is that it does not account for the relative importance of a reserve contribution, which we clarify at the example in Fig.~\ref{fig:exampleForUnfairReward}.
\begin{figure}[htb]
 		\vspace{-0.2\baselineskip}
        \centering
        \includegraphics[width=0.95\columnwidth]{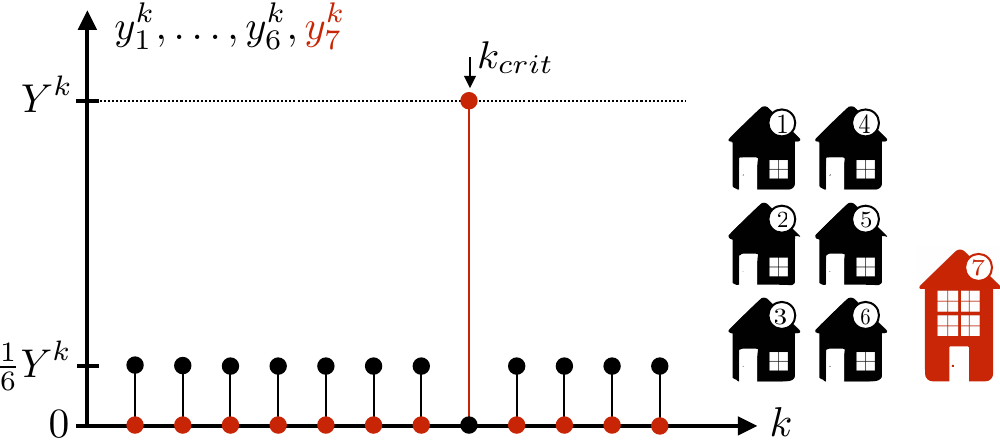}
        \vspace{-0.1cm}
        \caption{Illustrative example where seven buildings place an imbalanced reserve bid over a period of $N=12$ hours. During most of the time, six similar buildings (black) jointly offer reserves $y_1=y_2=\cdots = y_6$, while at \mbox{$k=k_{crit}$} a single building (red) determines the total reserve $Y^k=y_7^k$.}
        \label{fig:exampleForUnfairReward}
\end{figure}\vspace{-0.1\baselineskip}

In Fig.~\ref{fig:exampleForUnfairReward}, a single building determines the combined reserve offer at the critical time $k=k_{crit}$. We construct the example by using dynamic-free buildings and by constraining the black buildings' provision capability to zero at~$k_{crit}$.
Due to the time-constant constraint, the aggregation cannot provide any reserves if the red building leaves the group. 
Conversely, if a black building exits the aggregation, the reserve offer decreases only marginally (if at all), since the remaining black buildings can compensate for it.
Hence, due to its provision at~$k_{crit}$, the red building is more important. 
With the rewarding scheme~\eqref{eq:reservereward} however, the importance of provision during a particular time instance is not taken into account. 
The red building even earns less than any of the black buildings, e.g., for a time-constant price~$p$ it only gets $\nicefrac{6}{11}$ of each black building's profit.
We propose a Lagrangian-based approach to devise a reward allocation mechanism that better reflects the contribution of each building to the performance of the aggregation. As noted in Section~\ref{sec:aggregatedBiddingviaADMM}, the Lagrange multipliers $\lambda_b = (\lambda_b^1;\myDots;\lambda_b^{N})$ depend on time~$k$, building~$b$, and they also change with each ADMM iteration. Within Alg.~\ref{alg:ouradmm}, they orchestrate the buildings for achieving a constant combined reserve offer. 
 For a particular time instance, the intuition is that when the aggregation demand continually exceeds a building's proposition, i.e., $\bar y_b^k - y_b^k > 0$, then the Lagrangian update $$\lambda_b^k \leftarrow \lambda_b^k + \rho(\bar y_b^k - y_b^k)$$ accumulates~$\lambda_b^k$ towards a large value. 
 The more difficult the provision, the larger $\lambda_b^k$ grows before we reach a consensus. Fig.~\ref{fig:balance} illustrates the mechanism for $\lambda^k_b\geq0$.  
 In case of negative Lagrange multipliers, e.g., due to the initialization,  it can be seen from inspecting~\circled{1} and~\circled{2} in Alg.~\ref{alg:ouradmm} that  we recover a non-negative Lagrange multiplier quickly.

\begin{figure}[htb]
        \centering
        \includegraphics[width=1\columnwidth]{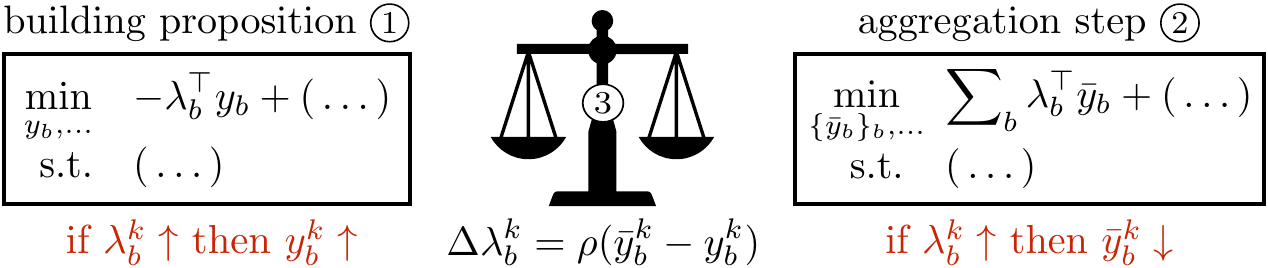}
        \vspace{-\baselineskip}
        \vspace{-0.05cm}
        \caption{Lagrangian mediation for $\lambda^k_b\geq0$. If $\bar y^k_b > y^k_b$, then $\lambda_b^k$ increases by $\Delta \lambda^k_b$ between ADMM iterations. This increase produces an incentive for raising $y_b^k$ and lowering $\bar y_b^k$, which reduces the gap between them. The amount of change in $y_b^k$ and $\bar y_b^k$ depends on the remaining constraints and objectives.}
        \label{fig:balance}
\end{figure}

According to the mechanism in Fig.~\ref{fig:balance}, the Lagrange multipliers indicate the time instances where the reserve provision is difficult. We use this effect to increase the reward during such critical hours. More specifically, we define the Lagrangian-based reward allocation
\begin{IEEEeqnarray}{rCl}
  \IEEEyesnumber \label{eq:Lagrangianreward} 
r_{b}^\Lambda \coloneqq \Lambda^\top y_b, 
\end{IEEEeqnarray}
where $\Lambda \in \mathbb R^N$ is an alternative price vector that replaces the static reward $p$ in~\eqref{eq:reservereward}. We describe the relation between $\Lambda$ and the Lagrange multipliers in the following proposition.
\begin{proposition}
\label{prop:LagrangianReward}
For Alg.~\ref{alg:ouradmm}, the following is true.
	\begin{itemize}
		\item[$(i)$] After each algorithm iteration, all buildings have the same Lagrange multiplier, i.e., $\lambda_1 = \dots = \lambda_M$. We write $\Lambda \coloneqq \lambda_b$ for any $b$, where $\Lambda=(\Lambda^1;\myDots; \Lambda^N)$. 
		\item[$(ii)$] The Lagrangian update (step \circled{3} in Alg.~\ref{alg:ouradmm}) can be written as $\Lambda \leftarrow \tfrac{\rho}{M}Y - \Omega$.
		\item[$(iii)$] For a converged solution $(\Lambda, Y, \{y_b\}_{b})$, the Lagrangian-based allocation \mbox{$r_{b}^\Lambda = \Lambda^\top y_b$} satisfies~$\sum\nolimits_{b} r^\Lambda_b = R$.
		\item[$(iv)$] For a feasible suboptimal solution $\left(Y^F,\{y_b^F\}_b\right)$ that is obtained with Prop.~\ref{prop:feasextract}, and
		\begin{IEEEeqnarray}{rCl}
  		\IEEEyesnumber \label{eq:aggregatedfeasible} \IEEEyessubnumber*
		\Omega^F &\coloneqq& \tfrac{\rho}{M} Y^F - \Lambda\\
		\Lambda^F &\coloneqq& \tfrac{1}{N}1_{N\times N}\left(\Omega^F+p\right)-\Omega^F,
		\end{IEEEeqnarray}
		where $\Lambda$ is obtained from the last algorithm iteration, the reward allocation $r_{b}^\Lambda = (\Lambda^F)^\top y_b^F$ satisfies $\sum\nolimits_{b} r^\Lambda_b =R$.
	\end{itemize}
\end{proposition}

 We show the proof of Prop.~\ref{prop:LagrangianReward} in Appendix~\ref{app:proofLagrangianReward}. 
In $(i)$, we define the alternative price vector $\Lambda$ as the building-invariant Lagrange multiplier. 
While the Lagrange multiplier is the same for all buildings, it still varies over the bidding time, indicating hours where provision is difficult. In Statement~$(ii)$, we provide a simple equation for~$\Lambda$. Statement $(iii)$ verifies that the allocation scheme distributes not more or less than the total reward~$R$, and Statement $(iv)$ makes it possible to use Lagrangian-based allocation after feasible~extraction. 
We interpret the building-invariance of the Lagrange multipliers as follows: $\lambda_b$ mediates  between building proposal~$y_b$ and aggregator demand~$\bar y_b$. However, the aggregator does not use any building-private information. In other words, it decides on the reserve demand~$\bar y_b$ without being able to differentiate between the buildings. Therefore, it always asks for the same change~$\bar y_b-y_b$, which results in the same multipliers.

The Lagrangian-based allocation provides an incentive to reduce provision imbalance. Hence, it strengthens the group's capability to provide reserves, which means that the aggregation grows more resilient against failures and configuration changes.
However, in some situations, a pure utilization of the new allocation scheme can be overly extreme for being economically acceptable, which becomes clear from Fig.~\ref{fig:lagrangianUnfairReward}. 
\begin{figure}[htb]
        \centering
        \vspace{-0.5\baselineskip}
        \includegraphics[width=0.77\columnwidth]{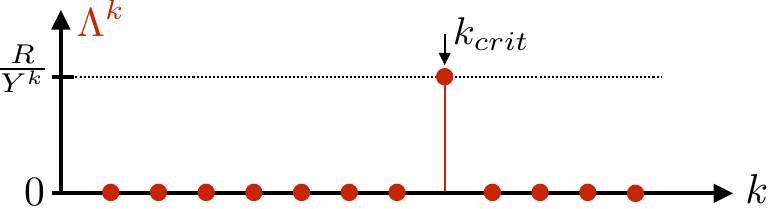}
        \vspace{-0.2cm}
        \caption{
        Lagrange multiplier $\Lambda$, which results from to the example in Fig.~\ref{fig:exampleForUnfairReward}. The reward concentrates in one point as a constraint limits the reserve provision only in a single time instant and the buildings are dynamic-free.
        }
        \label{fig:lagrangianUnfairReward}
\end{figure}

 While the extreme reward concentration in Fig.~\ref{fig:lagrangianUnfairReward} is specific to the constructed example, unrewarded hours can also occur in more generic setups. To alleviate this, we can use a mixed reward
	$\alpha\, r_b + (1-\alpha)\, r_b^\Lambda$
with an acceptable weight $\alpha \in [0,1]$.\looseness=-1


\section{Simulation Study}
\label{sec:simulation}
We apply ADMM-based reserve bidding in simulation, considering a $24$-hour bidding period with hourly discretization. We use a test set of $300$ building models, each generated by random perturbation from one of the following prototypes.\looseness=-1
\begin{itemize}
 	\item[$(i)$] Small building, obtained from~\cite[Section 4.5]{Oldewurtel_PhD2011}
 	\begin{itemize}
 		\item $3$ states for temperatures, using comfort constraints between $21^\circ$C and $25^\circ$C 
when the building is occupied	
\item $4$ constrained inputs for radiator, cooled ceiling, floor heating, and mechanical ventilation
 		\item $3$ modeled disturbances for outside temperature, solar radiation, and occupancy
 	\end{itemize}
 	\item[$(ii)$] Medium-sized building, obtained from~\cite{darivianakis2016power}
 	\begin{itemize}
 		\item $33$ states for temperatures, using comfort constraints between $20^\circ$C and $28^\circ$C  when the building is occupied
 		\item $5$ constrained inputs: $4$ for blinds and $1$ for heating
 		\item $7$ modeled disturbances: $1$ for occupancy, $2$ for ambient and ground temperature, and $4$ for solar radiation
    \end{itemize}
    \item[$(iii)$] Large building, obtained from~\cite{darivianakis2016power}
 	\begin{itemize}
 		\item $113$ states for temperatures, using comfort constraints between $20^\circ$C and $28^\circ$C when the building is occupied
 		\item $9$ constrained inputs: $4$ for blinds and $5$ for heating
 		\item $11$ modeled disturbances: $5$ for occupancy, $2$ for ambient and ground temperature, and $4$ for solar radiation
 	\end{itemize}
 \end{itemize}

 We classify the buildings in residential (mainly occupied during nighttime) and commercial (mainly occupied during daytime). Usually, the buildings provide most reserves while being unoccupied, as they are less restricted then.
 To obtain tractable bidding problems, we use linear decision rules as in 
Appendix~\ref{app:defCi} and\mbox{\cite[Section 4]{ZhangKamgGeorghiouGoulLyg_Aut17}}. We use free licenses for Yalmip~\cite{lofberg2004yalmip}, Gurobi~\cite{gurobi}, Julia~\cite{bezanson2017julia}, and JuMP~\cite{DunningHuchetteLubin2017}.\looseness=-1

\subsection{Example Scenario for Six Buildings}
We sample six buildings from the test set and obtain an aggregated reserve bid with Alg.~\ref{alg:ouradmm}. Fig.~\ref{fig:sixbuildingexample} illustrates the resulting bids~$y_1$ to~$y_6$ (bottom), and the Lagrange multiplier~$\Lambda$ (top).\looseness=-1 
\begin{figure}[h!tb]
		\vspace{-0.5\baselineskip}
        \centering
        \includegraphics[width=	.99\columnwidth]{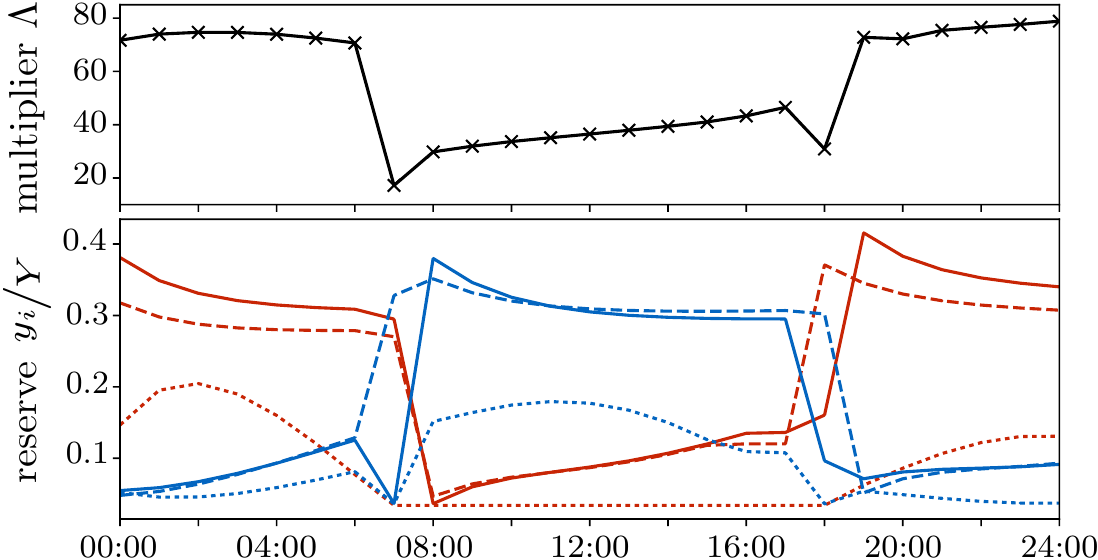}
        \vspace{-1.2\baselineskip}
        \caption{Example scenario of six buildings placing a joint reserve offer by solving~\eqref{eq:aggregatedProgram}. The horizontal axis shows the bidding period. The day-night provision pattern identifies three residential (blue) and three commercial (red) buildings.\looseness=-1}
        \label{fig:sixbuildingexample}
\end{figure}

 In Fig.~\ref{fig:sixbuildingexample}, we show the bids $y_{1},\myDots,y_6$ as a fraction of the combined bid~$Y$ (not illustrated), which adds up to be time-constant. In a separate simulation, we solve the individual problem~\eqref{eq:individualProgram} for each building, and we compare the result to the solution of the aggregated problem as illustrated in Fig.~\ref{fig:sixbuildingexample}. We observe that aggregated bidding with~\eqref{eq:aggregatedProgram} yields a combined bid $Y$ that is~$57\%$ larger than the accumulation of individual solutions from~\eqref{eq:individualProgram}.
This considerable aggregation advantage results from having a mixed group of residential and commercial buildings, which complement each other for satisfying the time-constant constraint. For the Lagrange multiplier, we see a day-night pattern that indicates a provision shortage at night. Hence, the commercial buildings are under-represented, and Lagrangian-based allocation places an incentive for more commercial buildings to join the group.\looseness=-1

\subsection{Analysis of the Necessary Number of ADMM Iterations}
\label{sec:stopppingCriteria}
The feasible extraction scheme presented in Section~\ref{sec:feasExtract} provides the opportunity to terminate the algorithm before it converges to a feasible solution. Therefore, the decision for algorithm termination only depends on the objective value 
\begin{IEEEeqnarray}{rCl}
  J\left( \{u_b(\cdot)\}_b, Y\right) \coloneqq \sum \limits_{b=1}^M c_b^\top  u_{b}(0) - p^\top  Y.
\end{IEEEeqnarray}
Fig.~\ref{fig:ItVsM} shows how many ADMM iterations are necessary until the feasible-extraction objective \mbox{$J^F \coloneqq J\left( \{u^F_b(\cdot)\}_b, Y^F\right)$} comes within $1\%$ of its precomputed optimal value~$J^\star>0$.
\begin{figure}[h!tb]
        \centering
        \includegraphics[width=.99\columnwidth]{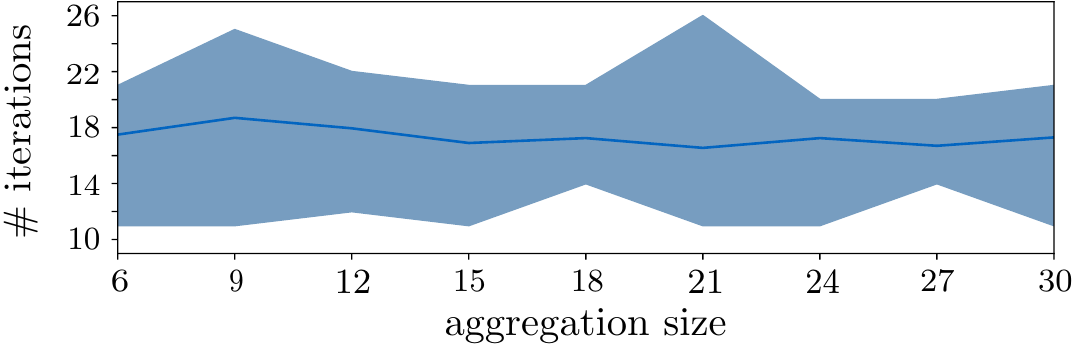}
        \vspace{-1.2\baselineskip}
        \caption{Required iterations of Alg.~\ref{alg:ouradmm} for reaching \mbox{${|J^F-J^\star|}/{J^\star} \leq 1\%$}  depending on the aggregation size $M$. For each aggregation size, we sample ten combinations of buildings from the building test set, and we show the average and range of the necessary number of ADMM iterations.}
        \label{fig:ItVsM}
        \iflongversion
        \else
        \vspace{-0.8\baselineskip}
        \fi
\end{figure}

 In Fig.~\ref{fig:ItVsM}, we observe: $(i)$ on average only $18$ ADMM iterations are needed to come within $1\%$ of the best-achievable objective; $(ii)$~the range between different instances of the same aggregation size is about $\pm 7$ iterations; and $(iii)$~the ADMM convergence characteristics (average and range) are fairly constant with the aggregation size. We expect the low number of iterations~$(i)$ and the small range~$(ii)$ as we have concentrated the main computational effort into the building-individual updates, which leaves a comparably simple and invariant task to the ADMM-based negotiation process. 
For the invariance to the aggregation size~$(iii)$, we note that the closed-form aggregation step~\eqref{eq:analyticaggregatorupdate} only depends on averaged information~$\Omega$. Hence, the mediation effort is not expected to vary with the number of buildings that contribute to this average.
With feasible extraction and the observations from Fig.~\ref{fig:ItVsM}, we can be confident to use the simplest possible ADMM termination criterion, namely just stopping after a fixed number of iterations (e.g.,~$25$). Besides being conceptually simple, this stopping criterion also leads to a fixed number of communication rounds between the buildings, making the entire decision process well-timed and easy to handle. In particular, given that we impose a time-limit for solving the building-individual problem, the total negotiation time for obtaining a joint reserve bid is constant.\looseness=-1

\section{Conclusion}
We discuss how flexible energy consumers, particularly buildings, can collaborate to provide frequency reserves in a reserve market setting. We show that when we use a distributed ADMM-based approach, we obtain conceptual benefits, such as low set-up costs, fully-decentralized computation, scalability, privacy, flexibility, and autonomy in the participating buildings.
We also show that ADMM provides a reward allocation scheme that reduces provision imbalance and therefore places an incentive to strengthen the group.
Finally, we emphasize the practicality of our approach by describing feasible extraction schemes, leading to a predictable stopping criterion and therefore a constant time for obtaining a joint reserve bid.\looseness=-1


%

\appendices

\section{Tractable Formulation of Building Constraints}\label{app:defCi}
In~\cite{ZhangKamgGoulartLyg_CDC2014, ZhangKamgGeorghiouGoulLyg_Aut17}, three steps are taken to obtain a tractable  reformulation of~\eqref{eq:consumerConstraintExample}. First, the uncertain reserve request \mbox{$s\in[-y,y]$} is represented through the normalized disturbance~$\zeta = (\zeta^1; \myDots ; \zeta^N)\in[-1,1]\subset\mathbb R^N$, such that~$s=s(\zeta) = \textnormal{diag}(y)\,\zeta$. This parametrization has the advantage that the normalized uncertainty set $[-1,1]$ is independent of the decision variable~$y$, which simplifies the maximization of the uncertainty set~$[-y,y]$\mbox{\cite[Section 3]{ZhangKamgGoulartLyg_CDC2014}}. In a second step, the control policies~$u(s)$ and~$\Delta u(s)$ are reparametrized with respect to the normalized disturbance as~$\tilde u(\zeta)$ and~$\Delta\tilde u(\zeta)$, which is still non-restrictive \cite[Prop.~3]{ZhangKamgGoulartLyg_CDC2014}. In a third (restrictive) step, affine parametrizations~$\tilde u(\zeta) = K\zeta + \kappa$ and~$\Delta\tilde u(\zeta) = F\zeta$ are used, where $K,F\in \mathcal L \subseteq \mathbb R^{Nm\times N}$, $\mathcal L$ is the set of lower block-diagonal matrices, and~$\kappa \in \mathbb R^{Nm}$ is a vector. The optimization is then performed over $K$, $F$, $\kappa$ and $y$, which results in
\begin{IEEEeqnarray}{l}
  \IEEEyesnumber \label{eq:consumerConstraintDetail} 
{\mathcal C} = \big\{\, \left( u(\cdot),y\right)  \,\,\big|\,\, \exists\, (K, F, \kappa, \tilde u(\cdot), \Delta  \tilde u(\cdot),  x(\cdot)) \,\, \text{ s.t.} \nonumber \\
 \tilde u(\zeta) = u(\textnormal{diag}(y)\,\zeta) = K \zeta + \kappa, \,\,\, \Delta \tilde u(\zeta) = F \zeta, \,\,\,  K,F\in\mathcal L, \nonumber\\
  \left. \begin{array}{r}
  x( \zeta) =  \bm A x^1 +  \bm B[ \tilde u(\zeta) +{\Delta   \tilde u(\zeta)}] +  \bm E v,\\
  \left( x(\zeta),  \tilde u(\zeta) + \Delta  \tilde u(\zeta)\right) \in \mathcal X\times {\mathcal U},\\
    y = (I_N \otimes \eta^\top) F\zeta
\end{array} \hspace{-1mm} \right\} \forall \zeta \in [-1, 1] \big\},\nonumber\\[-2.5ex]
\end{IEEEeqnarray}
which is consistent with~\eqref{eq:consumerConstraintExample} if~$\mathcal F$ contains linear functions only. 
Note that~\eqref{eq:consumerConstraintDetail} has an infinite number of constraints. By following the procedure in \cite[Section 4]{ZhangKamgGeorghiouGoulLyg_Aut17}, the resulting problem can be written as an exact dual reformulation with finitely many constraints. While the problem is convex and finite-dimensional, its solution remains computationally challenging. The reason is that the parametrization of the input~$\tilde u(\cdot)$ with the matrix~$K$ introduces a large number of decision variables. One way to reduce the computational complexity, at the cost of an increasingly  conservative solution, is to impose structure on~$K$. Popular strategies include restricting~$K$ to be block-diagonal, or even 
setting~$K=0$, which means that $\tilde u$ is preallocated and only $\Delta \tilde u$ depends on the actual reserve request.\looseness=-1

\section{Proofs}
\label{app:proofs}
\subsection{Proof of Proposition~\ref{prop:ouradmmconvergence}}\label{app:proofouradmmconvergence}
Problem~\eqref{eq:standardADMM} is equivalent to~\eqref{eq:aggregatedProgram}. Finite-dimensionality of~$\mathcal C$ makes it possible to use the convergence result in~\cite{boyd2011distributed}. We need the functions~$f,g$ to have a nonempty, closed, and convex epigraph, which is given by~\eqref{eq:splittingDefinition} and convexity of~$\mathcal C_b$. Further, due to $\mathcal C_b \neq \emptyset$, there is a feasible solution with $y_b\hspace{-0.3ex}=\hspace{-0.3ex}Y\hspace{-0.3ex}=\hspace{-0.3ex}0$.\looseness=-1
\hfill \IEEEQED

\iflongversion
\subsection{Proof of Proposition~\ref{prop:analyticaggregatorupdate}}\label{app:proofPropAnalytic}
 To simplify the notation, we use \mbox{$y \coloneqq (y_1; \myDots; y_b)$} and similarly $\bar y$ and $\lambda$. The aggregation step then becomes
\begin{IEEEeqnarray}{rCl}
\IEEEyesnumber \label{eq:rewrittenaggregatorupdate}\IEEEyessubnumber*
\min\nolimits_{(Y^k, \bar y)} & \quad & \lambda^\top \bar y + \tfrac{\rho}{2}\| \bar y -  y \|_2^2 - p^\top 1_{N\times 1}Y^k\\
\text{s.t.} 	&& 1_{N\times 1} Y^k = (1_{1\times M} \otimes I_N) \bar y, \label{eq:constraintInRewritten}
\end{IEEEeqnarray}
where $Y = 1_{N\times 1}{Y^k}$, i.e., $Y$ is parametrized by a single decision variable~$Y^k\in\mathbb R$, ensuring that \mbox{$Y^1 = Y^2 = \myDots = Y^N$}. To avoid confusion with the solution of Alg.~\ref{alg:ouradmm}, we denote the solution to~\eqref{eq:rewrittenaggregatorupdate} as $(\bm Y^k, \bar{\bm y})$. 
According to~\cite[Section 16.2]{Nocedal2000}, this solution satisfies the KKT conditions
 $0 = 1_{1\times N}(\bm \eta + p)$,
 $\bm{\bar{y}} = y + \tfrac{1}{\rho} \left( -\lambda - ( 1_{M\times 1} \otimes I_N) \bm \eta \right)$, and
 $1_{N\times 1} \bm Y^k = (1_{1\times M} \otimes I_N) \bm{\bar{y}}$,
where $\bm \eta$ is the optimal value of the Lagrange multiplier associated to~\eqref{eq:constraintInRewritten}. By using the Schur-complement method in~\cite{Nocedal2000}, we show that the KKT system is solved by
\begin{IEEEeqnarray}{rCl}
\IEEEyesnumber \label{eq:solKKT}\IEEEyessubnumber*
\bm Y^k &=& \tfrac{1}{\rho N} \left(1_{1\times NM} (\rho y - \lambda)  + M1_{1\times N}p \right) \label{eq:solYk}\\
\bm \eta  &=& \tfrac{1}{M}  \left( \left(1_{1\times M} \otimes I_N\right) \left(\rho y - \lambda\right) - \rho 1_{N\times 1} \bm Y^k \right)\\
\bm{\bar y} &=& \tfrac{1}{M}  (M - 1_{M\times M} \otimes I_N) (y - \tfrac{1}{\rho}\lambda) 
+ \tfrac{1}{M} I_{NM\times 1} \bm Y^k, \quad\quad\,\, \label{eq:solbary}
\end{IEEEeqnarray}
where we can write~\eqref{eq:solbary} component-wise. By using $\Omega = \tfrac{1}{M}  (1_{1\times M} \otimes I_N) (\rho y_b - \lambda_b)$ we conclude the proof.\looseness=-1
\hfill \IEEEQED
\fi

\subsection{Proof of Proposition~\ref{prop:feasextract}}\label{app:proofPropFeasible}
 $Y^F$ satisfies the time-constant constraint~\eqref{eq:Yconst}. From~\eqref{eq:createfeasy} it is $\sum_{b=1}^M {(y^k_{b})}^F = Y^F$, i.e.,~\eqref{eq:Ydef} is satisfied. It remains to show that~\eqref{eq:createfeasu} has a solution for all $b$. Due to ${\sum\nolimits_{b=1}^M y_{b}^k} \geq {{(Y^k)}^F}$, equation~\eqref{eq:createfeasy} ensures that ${(y_{b}^k)}^F \leq y_{b}^k$. Hence, by construction of $\mathcal C$ in~\eqref{eq:consumerConstraintExample},~\eqref{eq:createfeasu} is always feasible.
\hfill \IEEEQED

\subsection{Proof of Proposition~\ref{prop:LagrangianReward}}\label{app:proofLagrangianReward}
Inserting~\eqref{eq:ybresult} in step~\circled{3} of Alg.~\ref{alg:ouradmm} yields $\lambda_b = \tfrac{\rho}{M} Y - \Omega$ for any~$b$, which proves $(ii)$. $(i)$ follows as $Y$ and $\Omega$ are aggregated quantities. In $(iii)$, we consider solutions with $\bar y_b = y_b$, which we combine with~\eqref{eq:Yresult},~\eqref{eq:ybresult} to $\Lambda = \tfrac{1}{N}1_{N\times N}(\Omega+p) - \Omega$. $(iv)$ gives the same relation for~$\Lambda^F, \Omega^F$. In both cases, we multiply~$1_{1\times N}$ from the left and use $\tfrac{1}{N}1_{1\times N}1_{N\times N} = 1_{1\times N}$. The result is $1_{1\times N} \Lambda^{F} = 1_{1\times N} p$. $(iii)$ and $(iv)$ follow from writing~$1_{1\times N}$ as sum and multiplying~$Y^{F}$ from the right.\looseness=-1
\hfill \IEEEQED

\ifCLASSOPTIONcaptionsoff
  \newpage
\fi



%

\addtolength{\textheight}{-0 cm}   


\bibliographystyle{IEEEtran}  
\bibliography{library}

\begin{thebibliography}{10}
\providecommand{\url}[1]{#1}
\csname url@samestyle\endcsname
\providecommand{\newblock}{\relax}
\providecommand{\bibinfo}[2]{#2}
\providecommand{\BIBentrySTDinterwordspacing}{\spaceskip=0pt\relax}
\providecommand{\BIBentryALTinterwordstretchfactor}{4}
\providecommand{\BIBentryALTinterwordspacing}{\spaceskip=\fontdimen2\font plus
\BIBentryALTinterwordstretchfactor\fontdimen3\font minus
  \fontdimen4\font\relax}
\providecommand{\BIBforeignlanguage}[2]{{%
\expandafter\ifx\csname l@#1\endcsname\relax
\typeout{** WARNING: IEEEtran.bst: No hyphenation pattern has been}%
\typeout{** loaded for the language `#1'. Using the pattern for}%
\typeout{** the default language instead.}%
\else
\language=\csname l@#1\endcsname
\fi
#2}}
\providecommand{\BIBdecl}{\relax}
\BIBdecl

\bibitem{rebours2007survey}
Y.~G. Rebours, D.~S. Kirschen, M.~Trotignon, and S.~Rossignol, ``A survey of
  frequency and voltage control ancillary services - part i: Technical
  features,'' \emph{Power Systems, IEEE Transactions on}, vol.~22, no.~1, pp.
  350--357, 2007.

\bibitem{makarov2009operational}
Y.~V. Makarov, C.~Loutan, J.~Ma, and P.~de~Mello, ``Operational impacts of wind
  generation on {C}alifornia power systems,'' \emph{Power Systems, IEEE
  Transactions on}, vol.~24, no.~2, pp. 1039--1050, 2009.

\bibitem{callaway2011achieving}
D.~S. Callaway and I.~A. Hiskens, ``Achieving controllability of electric
  loads,'' \emph{Proceedings of the IEEE}, vol.~99, no.~1, pp. 184--199, Jan
  2011.

\bibitem{oldewurtel2013framework}
F.~Oldewurtel, T.~Borsche, M.~Bucher, P.~Fortenbacher, M.~G.~Vaya, T.~Haring,
  J.~Mathieu, O.~M{\'e}gel, E.~Vrettos, and G.~Andersson, ``A framework for and
  assessment of demand response and energy storage in power systems,'' in
  \emph{Bulk Power System Dynamics and Control, Security and Control of the
  Emerging Power Grid}.\hskip 1em plus 0.5em minus 0.4em\relax IEEE, 2013, pp.
  1--24.

\bibitem{Piette2012IntelligentBuilding}
M.~A. Piette, J.~Granderson, M.~Wetter, and S.~Kiliccote, ``Responsive and
  intelligent building information and control for low-energy and optimized
  grid integration,'' in \emph{ACEEE 2012 Summer Study on Energy Efficiency in
  Buildings}, Pacific Grove, CA, Oct 2012.

\bibitem{swissgridOkt2015}
\BIBentryALTinterwordspacing
\emph{Grundlagen Systemdienstleistungsprodukte}, Swissgrid AG, Feb 2017.
  [Online]. Available:
  \url{https://www.swissgrid.ch/dam/swissgrid/experts/ancillary_services/Dokumente/D170214_AS-Products_V9R2_de.pdf}
\BIBentrySTDinterwordspacing

\bibitem{chapman2016algorithmic}
A.~C. Chapman, G.~Verbi{\v{c}}, and D.~J. Hill, ``Algorithmic and strategic
  aspects to integrating demand-side aggregation and energy management
  methods,'' \emph{IEEE Transactions on Smart Grid}, vol.~7, no.~6, pp.
  2748--2760, 2016.

\bibitem{VrettosTPS2016}
E.~Vrettos, F.~Oldewurtel, and G.~Andersson, ``Robust energy-constrained
  frequency reserves from aggregations of commercial buildings,'' \emph{IEEE
  Transactions on Power Systems}, vol.~31, no.~6, pp. 4272--4285, 2016.

\bibitem{glowinski1975approximation}
R.~Glowinski and A.~Marroco, ``Sur l'approximation, par {\'e}l{\'e}ments finis
  d'ordre un, et la r{\'e}solution, par p{\'e}nalisation-dualit{\'e} d'une
  classe de probl{\`e}mes de dirichlet non lin{\'e}aires,'' \emph{Revue
  fran{\c{c}}aise d'automatique, informatique, recherche op{\'e}rationnelle.
  Analyse num{\'e}rique}, vol.~9, no.~2, pp. 41--76, 1975.

\bibitem{boyd2011distributed}
S.~Boyd, N.~Parikh, E.~Chu, B.~Peleato, and J.~Eckstein, ``Distributed
  optimization and statistical learning via the alternating direction method of
  multipliers,'' \emph{Foundations and Trends in Machine Learning}, vol.~3,
  no.~1, pp. 1--122, 2011.

\bibitem{diekerhof2017hierarchical}
M.~Diekerhof, F.~Peterssen, and A.~Monti, ``Hierarchical distributed robust
  optimization for demand response services,'' \emph{IEEE Transactions on Smart
  Grid}, 2017.

\bibitem{Vrettos14ifac}
E.~Vrettos, F.~Oldewurtel, F.~Zhu, and G.~Andersson, ``Robust provision of
  frequency reserves by office building aggregations,'' \emph{IFAC Proceedings
  Volumes}, vol.~47, no.~3, pp. 12\,068--12\,073, 2014.

\bibitem{Hao2013}
H.~Hao, A.~Kowli, Y.~Lin, P.~Barooah, and S.~Meyn, ``Ancillary service for the
  grid via control of commercial building {HVAC} systems,'' in \emph{American
  Control Conference}.\hskip 1em plus 0.5em minus 0.4em\relax Washington, USA:
  IEEE, 2013.

\bibitem{Maasoumy2014}
M.~Maasoumy, B.~Sanandaji, A.~Sangiovanni-Vincentelli, and K.~Poolla, ``Model
  predictive control of regulation services from commercial buildings to the
  smart grid,'' in \emph{American Control Conference}.\hskip 1em plus 0.5em
  minus 0.4em\relax Portland, Oregon, USA: IEEE, 2014.

\bibitem{VrettosAndersson_TSE2016}
E.~Vrettos and G.~Andersson, ``Scheduling and provision of secondary frequency
  reserves by aggregations of commercial buildings,'' \emph{IEEE Transactions
  on Sustainable Energy}, vol.~7, no.~2, pp. 850--864, 2016.

\bibitem{chen2013mpc}
C.~Chen, J.~Wang, Y.~Heo, and S.~Kishore, ``{MPC}-based appliance scheduling
  for residential building energy management controller,'' \emph{IEEE
  Transactions on Smart Grid}, vol.~4, no.~3, pp. 1401--1410, 2013.

\bibitem{vrettos2013combined}
E.~Vrettos and G.~Andersson, ``Combined load frequency control and active
  distribution network management with thermostatically controlled loads,'' in
  \emph{International Conference on Smart Grid Communications}.\hskip 1em plus
  0.5em minus 0.4em\relax IEEE, Oct 2013, pp. 247--252.

\bibitem{ZhangKamgGoulartLyg_CDC2014}
X.~Zhang, M.~Kamgarpour, P.~Goulart, and J.~Lygeros, ``Selling robustness
  margins: A framework for optimizing reserve capacities for linear systems,''
  in \emph{Conference on Decision and Control}.\hskip 1em plus 0.5em minus
  0.4em\relax IEEE, 2014, pp. 6419--6424.

\bibitem{BalandatOldewurtelChenTomlin2014}
M.~Balandat, F.~Oldewurtel, M.~Chen, and C.~Tomlin, ``Contract design for
  frequency regulation by aggregations of commercial buildings,'' in
  \emph{Allerton Conference on Communication, Control, and Computing}, Sep
  2014, pp. 38--45.

\bibitem{zhang2015stochastic}
X.~Zhang, E.~Vrettos, M.~Kamgarpour, G.~Andersson, and J.~Lygeros, ``Stochastic
  frequency reserve provision by chance-constrained control of commercial
  buildings,'' in \emph{European Control Conference}, Jul 2015, pp. 1134--1140.

\bibitem{GoreckiBitlisliogluStathopoulosJones_2015}
T.~T. Gorecki, A.~Bitlislioğlu, G.~Stathopoulos, and C.~N. Jones,
  ``Guaranteeing input tracking for constrained systems: Theory and application
  to demand response,'' in \emph{American Control Conference}.\hskip 1em plus
  0.5em minus 0.4em\relax IEEE, Jul 2015, pp. 232--237.

\bibitem{ZhangKamgGeorghiouGoulLyg_Aut17}
X.~Zhang, M.~Kamgarpour, A.~Georghiou, P.~Goulart, and J.~Lygeros, ``Robust
  optimal control with adjustable uncertainty sets,'' \emph{Automatica},
  vol.~75, pp. 249--259, 2017.

\bibitem{BitlisliogluGoreckiJones_TAC2017}
A.~Bitlislioglu, T.~T. Gorecki, and C.~N. Jones, ``Robust tracking
  commitment,'' \emph{IEEE Transactions on Automatic Control}, 2017.

\bibitem{bilgin2016provision}
E.~Bilgin, M.~C. Caramanis, I.~C. Paschalidis, and C.~G. Cassandras,
  ``Provision of regulation service by smart buildings,'' \emph{IEEE Trans. on
  Smart Grid}, vol.~7, no.~3, pp. 1683--1693, 2016.

\bibitem{taha2017buildings}
A.~F. Taha, N.~Gatsis, B.~Dong, A.~Pipri, and Z.~Li, ``Buildings-to-grid
  integration framework,'' \emph{IEEE Transactions on Smart Grid}, 2017.

\bibitem{kraning2014dynamic}
M.~Kraning, E.~Chu, J.~Lavaei, S.~Boyd \emph{et~al.}, ``Dynamic network energy
  management via proximal message passing,'' \emph{Foundations and Trends in
  Optimization}, vol.~1, no.~2, pp. 73--126, 2014.

\bibitem{weckx2015load}
S.~Weckx and J.~Driesen, ``Load balancing with ev chargers and pv inverters in
  unbalanced distribution grids,'' \emph{IEEE Transactions on Sustainable
  Energy}, vol.~6, no.~2, pp. 635--643, 2015.

\bibitem{hou2017distributed}
X.~Hou, Y.~Xiao, J.~Cai, J.~Hu, and J.~E. Braun, ``Distributed model predictive
  control via proximal jacobian {ADMM} for building control applications,'' in
  \emph{American Control Conference}.\hskip 1em plus 0.5em minus 0.4em\relax
  IEEE, 2017, pp. \mbox{37--43}.

\bibitem{verschae2016coordinated}
R.~Verschae, H.~Kawashima, T.~Kato, and T.~Matsuyama, ``Coordinated energy
  management for inter-community imbalance minimization,'' \emph{Renewable
  Energy}, vol.~87, pp. 922--935, 2016.

\bibitem{burger2017generation}
E.~M. Burger and S.~J. Moura, ``Generation following with thermostatically
  controlled loads via alternating direction method of multipliers sharing
  algorithm,'' \emph{Electric Power Systems Research}, vol. 146, pp. 141--160,
  2017.

\bibitem{OldewurtelenergyBuildings2012}
F.~Oldewurtel, A.~Parisio, C.~Jones, D.~Gyalistras, M.~Gwerder, V.~Stauch,
  B.~Lehmann, and M.~Morari, ``Use of model predictive control and weather
  forecasts for energy efficient building climate control,'' \emph{Energy and
  Buildings}, vol.~45, pp. 15--27, 2012.

\bibitem{Lehmann_2013}
B.~Lehmann, D.~Gyalistras, M.~Gwerder, K.~Wirth, and S.~Carl, ``Intermediate
  complexity model for model predictive control of integrated room
  automation,'' \emph{Energy and {Buildings}}, vol.~58, pp. 250--262, 2013.

\bibitem{Oldewurtel_PhD2011}
F.~Oldewurtel, ``Stochastic model predictive control for energy efficient
  building climate control,'' Ph.D. dissertation, ETH Zurich, 2011.

\bibitem{zhang:schildbach:sturzenegger:morari:13}
X.~Zhang, G.~Schildbach, D.~Sturzenegger, and M.~Morari, ``Scenario-based {MPC}
  for energy-efficient building climate control under weather and occupancy
  uncertainty,'' in \emph{European Control Conference}, Jul 2013, pp.
  1029--1034.

\bibitem{Darivianakis_CDC2015}
G.~Darivianakis, A.~Georghiou, R.~Smith, and J.~Lygeros, ``A stochastic
  optimization approach to cooperative building energy management via an energy
  hub,'' in \emph{Conference on Decision and Control}.\hskip 1em plus 0.5em
  minus 0.4em\relax Osaka, Japan: IEEE, Dec 2015, pp. 7814--7819.

\bibitem{Goulart2006}
P.~J. Goulart, E.~C. Kerrigan, and J.~M. Maciejowski, ``Optimization over state
  feedback policies for robust control with constraints,'' \emph{Automatica},
  vol.~42, no.~4, pp. 523--533, 2006.

\bibitem{darivianakis2016power}
G.~Darivianakis, A.~Georghiou, R.~Smith, and J.~Lygeros, ``The power of
  diversity: Data-driven robust predictive control for energy efficient
  buildings and districts,'' \emph{arXiv preprint arXiv:1607.05441}, 2016.

\bibitem{Scokaert1998}
P.~O. Scokaert and D.~Mayne, ``Min-max feedback model predictive control for
  constrained linear systems,'' \emph{IEEE Transactions on Automatic Control},
  vol.~43, no.~8, pp. 1136--1142, 1998.

\bibitem{LofbergPhD2003}
J.~L{\"o}fberg, \emph{Minimax approaches to robust model predictive
  control}.\hskip 1em plus 0.5em minus 0.4em\relax Link{\"o}ping University
  Electronic Press, 2003, vol. 812.

\bibitem{rey2016}
F.~Rey, P.~Hokayem, and J.~Lygeros, ``A tailored {ADMM} approach for power
  coordination in variable speed drives,'' in \emph{Preprints of the 20th IFAC
  World Congress}, Jul 2017, pp. 7664--7669.

\bibitem{lofberg2004yalmip}
J.~Lofberg, ``Yalmip: A toolbox for modeling and optimization in matlab,'' in
  \emph{Computer Aided Control Systems Design, 2004 IEEE International
  Symposium on}.\hskip 1em plus 0.5em minus 0.4em\relax IEEE, 2004, pp.
  284--289.

\bibitem{gurobi}
\BIBentryALTinterwordspacing
I.~Gurobi~Optimization, ``Gurobi optimizer reference manual,'' 2016. [Online].
  Available: \url{http://www.gurobi.com}
\BIBentrySTDinterwordspacing

\bibitem{bezanson2017julia}
J.~Bezanson, A.~Edelman, S.~Karpinski, and V.~B. Shah, ``Julia: A fresh
  approach to numerical computing,'' \emph{SIAM Review}, vol.~59, no.~1, pp.
  65--98, 2017.

\bibitem{DunningHuchetteLubin2017}
I.~Dunning, J.~Huchette, and M.~Lubin, ``{JuMP}: A modeling language for math.
  optimization,'' \emph{SIAM Review}, vol.~59, no.~2, pp. 295--320, 2017.

\bibitem{Nocedal2000}
J.~Nocedal and S.~J. Wright, \emph{Numerical optimization}.\hskip 1em plus
  0.5em minus 0.4em\relax Springer Science \& Business Media, 1975, vol.~9,
  no.~4.

\end{thebibliography}

\vspace{-4.5\baselineskip}

\vfill

\begin{IEEEbiography}[{\includegraphics[width=1in,height=1.25in,clip,keepaspectratio]{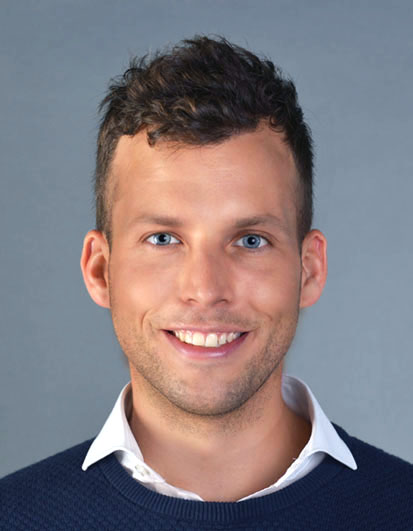}}]{Felix Rey} is a Ph.D. candidate in the Automatic Control Laboratory at ETH Zurich.
He received his B.Eng. degree in Electrical Engineering and Information Technology from the University of Applied Sciences Constance (Germany, 2011), where he became a fellow of the German National Academic Foundation. He obtained his M.Sc. degree in Electrical Engineering and Information Technology from the Karlsruhe Insitute of Technology (Germany, 2014). In 2014 he joined the Automatic Control Laboratory at ETH Zurich as doctoral student. His research interests include model predictive control, as well as distributed and embedded optimization, in particular with ADMM.\looseness=-1
\end{IEEEbiography}

\begin{IEEEbiography}[{\includegraphics[width=1in,height=1.25in,clip,keepaspectratio]{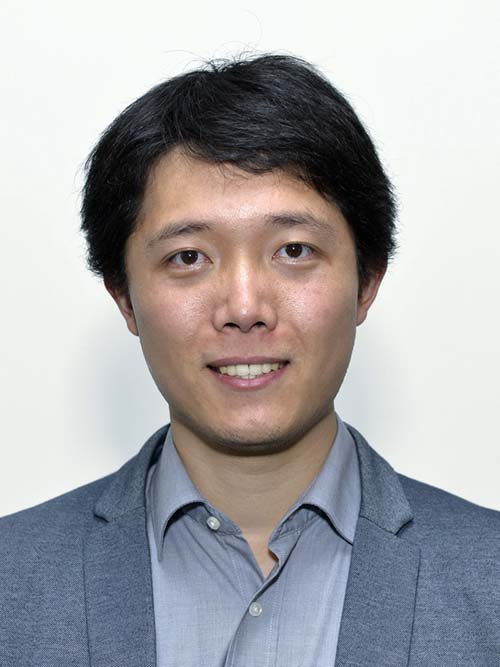}}]{Xiaojing (George) Zhang} 
received his B.Sc and M.Sc. degrees in electrical engineering and information technology, in 2010 and 2012, respectively, from the Swiss Federal Institute of Technology (ETH) Zurich, Switzerland. In 2016, he received his Ph.D. degree in automatic control, also from ETH Zurich. He is currently a post-doctoral fellow and the Associate Director at the Hyundai Center of Excellence at the University of California, Berkeley. His research interests include modeling, analysis, and control of stochastic uncertain systems, randomized and data-driven algorithms, as well as stochastic optimization, with application to autonomous vehicles, energy-efficient buildings and electric power systems.
\end{IEEEbiography}
\vspace{-1.92\baselineskip}
\begin{IEEEbiography}[{\includegraphics[width=.95in,height=1.25in,clip,keepaspectratio]{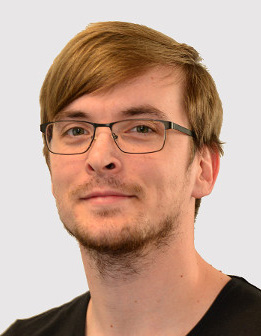}}]{Sandro Merkli}  received his BSc. and MSc. degrees in electrical engineering in 2012 and 2014, respectively, from ETH Zurich. He is currently a Ph.D. candidate at the Automatic Control Laboratory at ETH Zurich. His research interests are numerical optimization in general as well as its applications in power system operation and design. 
\end{IEEEbiography}
\vspace{-1.92\baselineskip}
\begin{IEEEbiography}[{\includegraphics[width=1in,height=1.25in,clip,keepaspectratio]{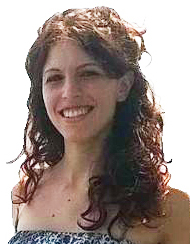}}]{Valentina Agliati} received her BSc. and MSc. degrees in automation and control engineering in 2014 and 2017, respectively, from Politecnico di Milano. In 2016 she wrote her master's thesis in the Automatic Control Laboratory at ETH Zurich as an exchange student.
\end{IEEEbiography}
\vspace{-1.92\baselineskip}
\begin{IEEEbiography}[{\includegraphics[width=.95in,height=1.25in,clip,keepaspectratio]{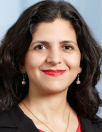}}]{Maryam Kamgarpour}
 is an Assistant Professor in the Automatic Control Lab at ETH Zurich. She obtained her Master’s and Ph.D. in Control Systems at the University of California, Berkeley (2007, 2011) and her Bachelor of Applied Sciences in Systems Design Engineering from University of Waterloo, Canada (2005). Her research is on safety verification and optimal control of large-scale uncertain dynamical systems with applications in air traffic and power grid systems. She is the recipient of the NASA High Potential Individual Award, the NASA Excellence in Publication Award (2010) and the European Union (ERC) Starting Grant 2015.
\end{IEEEbiography}
\vspace{-1.92\baselineskip}
\begin{IEEEbiography}[{\includegraphics[width=1in,height=1.25in,clip,keepaspectratio]{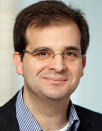}}]{John Lygeros}
  received the B.Eng. degree in electrical engineering in 1990 and an M.Sc. degree in Systems and Control in 1991, both at Imperial College of Science Technology and Medicine, London, UK, in 1990 and 1991, respectively, and the Ph.D. degree from the Electrical Engineering and Computer Sciences Department, University of California, Berkeley, in 1996. He holds the chair of Computation and Control at the Swiss Federal Institute of Technology (ETH) Zurich, Switzerland, where he is currently serving as the Head of the Automatic Control Laboratory. During the period 1996–-2000 he held a series of research appointments at the National Automated Highway Systems Consortium, Berkeley, the Laboratory for Computer Science, M.I.T., and the Electrical Engineering and Computer Sciences Department at U.C. Berkeley. Between 2000 and 2003 he was a University Lecturer at the Department of Engineering, University of Cambridge, UK, and a Fellow of Churchill College. Between 2003 and 2006 he was an Assistant Professor at the Department of Electrical and Computer Engineering, University of Patras, Greece. In July 2006 he joined the Automatic Control Laboratory at ETH Zurich, first as an Associate Professor, and since January 2010 as Full Professor. His research interests include modeling, analysis, and control of hierarchical, hybrid, and stochastic systems, with applications to biochemical networks, automated highway systems, air traffic management, power grids and camera networks. John Lygeros is a Fellow of the IEEE and a Member of the IET and of the Technical Chamber of Greece.\looseness=-1
\end{IEEEbiography}




\end{document}